\documentclass[11pt]{amsart}
\usepackage{cite}
\usepackage{wrapfig, lipsum,booktabs}
\usepackage{graphicx}
\usepackage{amssymb}
\usepackage{epstopdf}
\usepackage{verbatim}
\usepackage{bm}
\usepackage{multicol}
\usepackage{multirow}
\usepackage{subfigure}
\usepackage{fancyhdr} 
\usepackage{float}
\usepackage{stmaryrd}
\usepackage{color}
\usepackage{soul}
\usepackage{tikz}
\usetikzlibrary{patterns}
\usepackage{pgfplots}
\usepackage{cancel}

\usepackage{multirow}
\usepackage{amsmath,amssymb,eucal}
\usepackage{graphicx,subfigure,epsfig}
\usepackage{psfrag}
\usepackage{url}
\usepackage{todonotes}
\usepackage[top=1in, bottom=1.in, left=1in, right=1in]{geometry}

\begin{document}
\title[div-free HDG for CHNS]{A divergence-free HDG scheme for 
the Cahn-Hilliard phase-field model for two-phase incompressible flow}
\author{Guosheng Fu}
\address{Department of Applied and Computational Mathematics and 
Statistics, University of Notre Dame, USA.}
\email{gfu@nd.edu}
% \thanks{The author gratefully acknowledges the partial support of this work
% under AFOSR contract FA9550-12-1-0399.}

\begin{abstract}
  We construct a divergence-free HDG scheme
   for the Cahn-Hilliard-Navier-Stokes phase field model. 
The scheme is robust in the convection-dominated regime, produce a globally
divergence-free velocity approximation, and can be efficiently 
implemented via static condensation.
Two numerical benchmark problems, namely the rising bubble problem, and the 
Rayleigh-Taylor instability problem are used to show the good performance
of the proposed scheme.
\end{abstract}
\maketitle

\section{Introduction}
\label{sec:intro}
In the past few decades there has been tremendous progress in the development and analysis of numerical methods for single-phase incompressible flow problems. Recently, the importance of producing an exactly divergence-free velocity approximation for single-phase incompressible Stokes and Navier-Stokes equations 
within the mixed finite element framework has been stressed \cite{JohnLinkeMerdonNeilanRebholz17}.
There are active research on the extension of numerical methods developed for single-phase incompressible flow problems to two-phase incompressible flow problems \cite{GrossReusken08}. 
The fundamental issue relevant for the simulation of two-phase flows that is non-existent in one-phase incompressible flow problems is the numerical treatment of the {\it unknown} interface. 

%%% NEW COMMANDS
\newcommand{\pol}{\mathcal{P}} 
\newcommand{\polq}{\mathcal{Q}} 
%
% boldsymbols
\newcommand{\bld}[1]{\hbox{\boldmath$#1$}}    
\newcommand{\Th}{\mathcal{T}_h}
\newcommand{\Eh}{\mathcal{E}_h}
\newcommand{\jmp}[1]{[\![#1 ]\!]}
\newcommand{\Oh}{\mathcal{T}_h}
\newcommand{\gfnote}[1]{\todo[inline,color=orange!20]{GF: #1}}

In this paper, we initalize an investigation on the application  of 
the divergence-free hybridizable discontinuous Galerkin (HDG) methods \cite{Lehrenfeld10,LehrenfeldSchoberl16} 
to  a diffusive interface model based on the Cahn-Hilliard equations  for the two-phase incompressible flow. 
In particular, we consider the following Navier-Stokes-Cahn-Hilliard system
 proposed in \cite{DingSpeltShu07}:
\begin{subequations}
  \label{ns-eq}
  \begin{alignat}{2}
    \label{ns-eq1}
    \partial_t{\phi} +\bld u\cdot \nabla \phi =&\; \nabla\cdot
    (M(\phi)\nabla\mu),\quad&&\text{ in } \Omega,\\
    \mu =&\; \widetilde{\sigma}(\epsilon^{-1}W'(\phi)-\epsilon\triangle \phi)
    %+\frac{1}{2}\rho'(\phi)|\bld u|^2
    ,\quad&&\text{ in } \Omega,\\
%    \sqrt{\rho}\partial_t(\sqrt{\rho}\bld u)
    \label{ns-eq2}
    \rho(\phi)(\partial_t(\bld u) 
    +\nabla\cdot(\bld u\otimes\bld u))
%    -
  %  \frac12\nabla\cdot(\rho\bld u)\bld u
    =&\;
    \nabla \cdot
%    \widetilde
    %{\bld \sigma}
    (2\nu(\phi)\mathbf{D}(\bld u))-\nabla {p}
    +\rho\bld f+\mu\nabla \phi
    , \quad&& \text{ in } \;\Omega,\\
    \label{ns-div}
  \nabla \cdot\bld u = &\;0,\quad     && \text{ in
    } \;\Omega,\\
    \frac{\partial \phi}{\partial \bld n}=
    \frac{\partial \mu}{\partial \bld n}=&\;0, \;\;
    \bld u = \bld 0,\quad     && \text{ on } \;\partial\Omega,
    %\quad in \;\;\Omega,
  \end{alignat}
\end{subequations}
where $\Omega\subset \mathbb{R}^2$ is a 
two-dimentional domain such that 
$\overline{\Omega}=\overline{\Omega}_1\cup\overline{\Omega}_2$
and $\Omega_1\cap\Omega_2=\emptyset$,
where $\Omega_i$ denotes the  subdomain
occupies an incompressible fluid with densitiy 
$\rho_i$ and dynamic viscosity $\nu_i$ for $i=1,2$.
Here $\mathbf{D}(\bld u):=\frac12(\nabla \bld u + (\nabla \bld u)^T)
$ is the symmetric small strain tensor,
%the stress tensor
%$\bld{\sigma}:=2\nu(\phi)\mathbf{D}(\bld u)
%-{p}\bld I$ .
and $\bld u, {p}, \phi, \mu$ are the velocity, pressure, phase field
variable, and chemical potential, respectively.
Furthermore, 
the variable density $\rho$ and viscosity $\nu$ are
slave variables of $\phi$ given by the linear average
\begin{align}
  \label{den-vis}
  \rho(\phi) := \rho_1(1+\phi)/2+\rho_2(1-\phi)/2,\quad
  \nu(\phi) := \nu_1(1+\phi)/2+\nu_2(1-\phi)/2.
\end{align}
%and the initial condition 
%for the phase field variable $\phi$ 
%is taken to be $1$ on the exterior domain $\Omega_1(0)$
%and $-1$ on the interior domain $\Omega_2(0)$.
%Hence, we have 
% \[
%   \rho'(\phi) = \frac{\rho_1-\rho_2}{2},
%\] 
%which is independent of the variable $\phi$.
The function $W(\phi)$ is a double well potential, 
here we use
\begin{align}
  \label{db-well}
  W(\phi):=\frac14(\phi^2-1)^2.
\end{align}
The function $M(\phi)$ is a mobility function,  here we use  
\begin{align}
  \label{mob}
  M(\phi) := \gamma (\phi^2-1)^2,
\end{align}
where $\gamma$ is a constant mobility coefficient.
Finally, 
$\epsilon$ defines a length scale over which the interface is smeared out,
and the parameter $\widetilde{\sigma}$ is a scaled surface tension, which
is related to the physical surface tension $\sigma$ by 
$ \widetilde{\sigma}:=\frac{3}{2\sqrt{2}}\sigma$.

It is noteworthy to straighten that the model 
\eqref{ns-eq}, unfortunately, does not admit an energy law.
Following \cite{GuermondQuartapelle00}, an energy estimate of the model
could be shown by introducing the new variable $r=\sqrt{\rho}$, see
\cite{ShenYang10}. 
We also mention that 
thermodynamically consistent phase-field models that admit an energy law
was recently derived in
\cite{Abels12}.
It was numerically shown in \cite{AlanVoigt12} that very similar results 
were obtained for standard finite element
discretizations for the model \eqref{ns-eq}  and the thermodynamically consistent  model in
\cite{Abels12} in the context of a rising bubble benchmark problem
\cite{Hysing09}. For this reason, we focus on the construction of HDG
scheme for the slightly simplier model \eqref{ns-eq}, although it does not
admit an energy law.

The rest of the paper is organized as follows.
In Section \ref{sec:scheme}, we first introduce 
the divergence-free HDG-based spatial discretization for the model problem 
\eqref{ns-eq}, then apply a standard Crank-Nilson based IMEX time
discretization for the resulting ODE system, which leads to a 
(conditionally stable) 
linear decoupled fully discrete scheme.
Then in Section \ref{sec:num}, we first numerically study the convergence
property of the proposed scheme, which indicates second order convergence
in time, and optimal $(k+1)$-th order of convergence
in space when polynomials of degree $k$
is used for the spatial discretization. Next, we apply our scheme to two
benchmark tests, namely, the rising bubble problem, and the Rayleigh–Taylor 
instability problem.
We conclude in Section \ref{sec:conclude} with some future work.

\section{The divergence-free HDG scheme}
\label{sec:scheme}
In this section, we introduce the divergence-free HDG scheme for the model
problem \eqref{ns-eq} in two dimensions. 
Although our  scheme can be defined on hybrid triangular/quadrilateral
meshes, in this paper we only present the scheme on structured rectangular
meshes. 
To this end, let $\Th:=\{T\}$ be a conforming 
rectangular triangulation of the rectangular domain $\Omega$, 
and let $\Eh=\{F\}$ be the collection
of edges of
$\Th$.
We set $h$ to be the maximum mesh size of $\Th$.
Given a rectangular element $T$, 
we denote  $\mathcal{P}^{m,n}(T)$ as 
the space of polynomials of degree at most $m$ in the first argument and 
at most $n$ in the second argument on the element $T$. 
We further denote $\mathcal{Q}^m(T):=\mathcal{P}^{m,m}(T)$ to simplify
notation.
Given an edge $F$,  
we denote  $\mathcal{P}^{m}(F)$ as 
the space of polynomials of degree at most $m$ on the edge $F$.
On each element $T$, 
we denote the tangential component of a vector $\bld v$ on an edge $F$
by $\mathsf{tang}(\bld v) := \bld v-(\bld v\cdot \bld n)\bld n$, where $\bld n$ 
is the unit normal vector on $F$.

The following finite element spaces will be used:
\begin{subequations}
  \label{space}
  \begin{align}
    \label{space-1}
    W^{k}_h : =&\;
    \{w\in L^2(\Th):\quad w\in\polq^k(T),\;\;\forall T\in\Th\}, \\
    %\prod_{T\in\Oh}\pol^{k}(T),\\
  X^{k}_h : =&\; 
  %\prod_{F\in\Eh}\pol^{k}(F)
  \{\widehat{w}\in H^1(\Eh):\quad \widehat{w}\in\pol^k(F),\;\;
  \forall F\in\Eh\}, \\
    \Phi^{k}_h : =&\;
    \{\xi\in H^1_0(\Th):\quad \xi \in\polq^k(T),\;\;\forall T\in\Th\}, \\
    \label{space-div}
    \bld V^k_h : =&\; \nabla \times \Phi^{k+1}_h,\\ 
  %\subset H_0(\mathrm{div}_0,\Omega),\\
    M^k_{h} := &\;\{\widehat{\bld v}\in [L^2(\Eh)]^2:
      \quad
      \widehat{\bld v}\in [\pol^{k}(F)]^2, 
        \widehat{\bld v}\cdot\bld n = 0 \,\forall F\in\Eh, \;\;
      \widehat{\bld v} = 0 \,\;\;\forall F\subset \partial\Omega\},
        \end{align}
      \end{subequations}
      where 
      $k$ is
      the polynomial degree, and the {\it curl} operator in the velocity 
      space \eqref{space-div} is the rotated gradient: 
      $\nabla \times = (\partial_y, -\partial_x).$
      Note that functions in $X^k_h$ are defined only on the mesh skeleton
      and are continuous on the mesh nodes, while functions in
      $M^k_h$ are defined only on the mesh skeleton 
      and have normal component {\it zero}.
      Note also that functions in the velocity space \eqref{space-div} is globally
      divergence-free, and $\{\nabla \times \xi_j\}_{j=1}^N$ forms a set of
      basis for $\bld V_h^k$, whenever $\{\xi_j\}_{j=1}^N$ forms a set of basis
      for the scalar continuous finite element space $\Phi_h^{k+1}$.
      In particular, integration-by-parts yields the following identity:
      \begin{align}
        \label{pressure}
        (\nabla p, \bld v) = -(p, \nabla \cdot\bld v) = 0,\quad 
        \forall v\in \bld V_h^k,
      \end{align}
      which will help us to exclude the presure approximation in the numerical
      scheme.

      The proposed spatial discretization is given as follows:
      Find $(\phi_h, \widehat{\phi}_h, 
      \mu_h,\widehat{\mu}_h, \bld u_h, \widehat{\bld u}_h)\in 
      W^k_h\times X^k_h\times W^k_h\times X^k_h\times \bld V^k_h\times
      M^{k-1}_h
     $
     such that 
      \begin{subequations}
        \label{semi}
        \begin{alignat}{2}
        \label{semi-phi}
          (\partial_t\phi_h, \psi)+
        \mathcal{C}_h^1\left((\bld u_h, \phi_h,\widehat{\phi}_h), 
        (\psi,\widehat{\psi})\right)
          +\mathcal{D}_h\left((M_h,\mu_h,\widehat{\mu}_h), (\psi,
          \widehat{\psi})\right)
          =&\;0,\\
        \label{semi-mu}
          -(\mu_h, \eta)+
          \widetilde{\sigma}\epsilon
          (W'(\phi_h), \eta) +
          \mathcal{D}_h\left((\widetilde{\sigma}\epsilon, 
            \phi_h, \widehat{\phi}_h), (\eta,
            \widehat{\eta})\right)
          =&\;0,\\
          \label{semi-u}
          \left(\rho_h\partial_t\bld u_h, \bld v\right)
          +\mathcal{C}_h^2\left((\rho_h, \bld u_h, \bld u_h),
          (\bld v,\widehat{\bld v})\right)
          %\quad\quad&\\
            +\mathcal{B}_h\left((\nu_h, \bld u_h, \widehat{\bld u}_h), (\bld v, \widehat{\bld v})
            \right)\quad&\;\\
            - \mathcal{C}_h^3\left((\phi_h,\widehat{\phi}_h,\mu_h), \bld v
            \right)
          = &\;
          (\rho_h\bld f, \bld v),\nonumber
        \end{alignat}
      for all 
      $(\psi, \widehat{\psi},\eta,\widehat{\eta},
      \bld v, \widehat{\bld v})\in 
      W^k_h\times X^k_h\times W^k_h\times X^k_h\times \bld V^k_h\times
      M^{k-1}_h
     $ where we denote 
      $\rho_h:=\rho(\phi_h), M_h:=M(\phi_h)$, and $\nu_h:=\nu(\phi_h)$ 
      to shorten the notation.
      Here $(\cdot,\cdot)$ denotes the $L^2$-inner product on the mesh $\Th$,
      and the operators in \eqref{semi} are given as follows:
      \begin{alignat*}{2}
        \mathcal{C}_h^1\left((\bld u_h, \phi_h,\widehat{\phi}_h),
        (\psi,\widehat{\psi})\right)
        :=& \sum_{T\in\Th}\left(-\int_T\bld u_h\phi_h\cdot\nabla\psi\, \mathrm{dx}-
    \int_{\partial T}\bld u_h\cdot\bld
    n\widehat{\phi}_h^{up}(\psi-\widehat{\psi})\mathrm{ds}
  \right),
  \\
 \mathcal{D}_h\left((c,\mu_h, \widehat{\mu}_h),
 (\psi,\widehat{\psi})\right)
          := &\sum_{T\in\Th}\Big(\int_Tc\nabla\mu_h\cdot\nabla\psi\, \mathrm{dx}
            -
   \int_{\partial T}c
   \frac{\partial \mu_h}{\partial n}
   (\psi-\widehat{\psi})
   \mathrm{ds}
  \\
             & \hspace{-.5cm}-
   \int_{\partial T}c
   \frac{\partial \psi}{\partial n}
   (\mu_h-\widehat{\mu}_h)
   \mathrm{ds}
   +
   \int_{\partial T}c
   \frac{\alpha(k+1)^2}{h} 
   (\mu_h-\widehat{\mu}_h)
   (\psi-\widehat{\psi})
   \,
   \mathrm{ds}
  \Big),\nonumber
  \\
 \mathcal{C}_h^2\left((\rho_h, \bld w_h,\bld u_h),
 \bld v)\right)
  :=& \sum_{T\in\Th}\Big(\int_T
    \rho_h\nabla\cdot(\bld w_h\otimes \bld u_h)\cdot \bld v
    \, \mathrm{dx} +
               \int_{\partial T}{\rho}_h\bld w_h\cdot\bld n
               (\widehat{\bld u}_h^{up}-\bld u_h)\cdot\bld v\,
   \mathrm{ds}
  \Big),\nonumber\\
  \mathcal{B}_h((c, \bld u_h, \widehat{\bld u}_h), 
        (\bld v, \widehat{\bld v}))
        :=&\; 
        2 \sum_{T\in\Th}\left(\int_T c%(\phi_h)
          \mathbf{D}(\bld u_h):
          \mathbf{D}(\bld v)
          \, \mathrm{dx}
          -\int_{\partial T}
          c
          \mathbf{D}(\bld u_h)\bld n \cdot 
          \mathsf{tang}(
            \bld v-\widehat{\bld v}
          %\jmp{\bld v}
          )
          \,\mathrm{ds}
        \right.\\
          &\;\hspace{-4.cm}
          \left.
            - \int_{\partial T}
            c
          \mathbf{D}(\bld v)\bld n \cdot 
            \mathsf{tang}(
            \bld u_h-\widehat{\bld u}_h 
            %\jmp{\bld u_h}
            )\,\mathrm{ds}
            + \int_{\partial T}
            c
            \frac{\alpha(k+1)^2}{h}
            \mathsf{tang}(\Pi_{k-1}
            \bld u_h-\widehat{\bld u}_h
            %\jmp{\bld u_h}
            )
            \cdot\mathsf{tang}(\Pi_{k-1}
            \bld v-\widehat{\bld v}
            %\jmp{\bld v}
            )
            \,\mathrm{ds}
          \right),\\
          \mathcal{C}_h^3\left((\phi_h, \widehat{\phi}_h
            ,\mu_h),
        \bld v)\right)
        :=& \sum_{T\in\Th}\left(-\int_T\phi_h\nabla\mu_h\cdot\bld v\,
          \mathrm{dx}+
    \int_{\partial T}\bld v\cdot\bld
    n\widehat{\phi}_h\mu_h\mathrm{ds}
  \right),
\end{alignat*}
      \end{subequations}
      where
      the {\it upwinding} numerical flux
      $\widehat{\bld u}_h^{up}|_F=\bld u_h^-|_F:= (\bld u_h|_{T^-})|_{F}$ where 
      $T^-$ is the  element sharing the edge $F$ such that 
      $(\bld u_h\cdot\bld n|_{T^-})|_{F}>0$, and the {\it upwinding}
      numerical flux
      \begin{align*}
        \widehat{\phi}_h^{up}:=\left\{
          \begin{tabular}{l}
            $\phi_h$, if $\bld u_h\cdot\bld n\ge 0$,\\[1.5ex]
            $\widehat{\phi}_h$, if $\bld u_h\cdot\bld n< 0$,
          \end{tabular}
        \right.
      \end{align*}
    Here  $\Pi_{k-1}$ is the $L^2$-projection operator onto the space 
    $M_h^{k-1}$, and 
      $\alpha>0$ is a sufficiently large stabilization constant
      to ensure positivity of the viscous operators 
    $\mathcal{D}_h((c,\cdot,\cdot), (\cdot,\cdot))$ and 
    $\mathcal{B}_h((c,\cdot,\cdot), (\cdot,\cdot))$,  
      which is taken to be $\alpha=4$ in all our numerical simulations.

      We remark that in the above scheme \eqref{semi}, 
      the embeded discontinuous Galerkin (EDG) approach 
\cite{CockburnGuzmanSoonStolarski09,NguyenPeraireCockburn15} is used 
to discretize the second-order terms in equations \eqref{semi-phi} and 
\eqref{semi-mu}. The convection term in \eqref{semi-phi} is
discretized using an EDG-based upwinding.
The convection term in the momentum equation \eqref{semi-u} is discretized
using a classical DG-based upwinding, the viscous term therein is
discretized using a divergence-free HDG approach with projected jumps
\cite{Lehrenfeld10,LehrenfeldSchoberl16} to further save computational
cost, the surface tension force term is obtained by integration by
parts,
%and the fact that $\nabla \cdot\bld v=0$ for all test functions 
%$\bld v\in \bld V_h^k$, 
and the pressure term cancels out due to the
equality \eqref{pressure}.
In particular, we emphasis that: 
\begin{itemize}
  \item [(i)]
the numerical dissipation (from upwinding)
in the treatment
of the two convection terms in \eqref{semi-phi} and \eqref{semi-u} is
beneficial to stabilize the scheme in the convection-dominated regime
without using any extra residual-based stabilization;
\item [(ii)]
due to 
the choice of the (globally)  divergence-free velocity space $\bld V_h^k$,
the 
divergence-free condition \eqref{ns-div} is always satisfied strongly for the 
scheme \eqref{semi}, which is excepted to be more robust than 
methods that only 
satisfy the divergence-free condition in an approximate/discrete sense,
c.f. \cite{JohnLinkeMerdonNeilanRebholz17}.
\item [(iii)]
the introduction of the hybrid unknonws 
(in the context of EDG for phase-field variables, and 
HDG with projected jumps for the velocity variable) is a natural way to
further improve the efficiency of the resulting linear system solvers 
(via static condensation).
\end{itemize}
%which 
%is the major reason for us to consider such DG-based scheme.

For the temporal discretization, we simply use the
second-order Crank-Nilson-Adams-Bashforth IMEX
approach \cite{AscherRuuthWetton93} as follows.
For any positive integer $j\in \mathbb{Z}_+$,
let $(\rho_h^{j-1}, \widehat{\rho}_h^{j-1}, \bld u_h^{j-1}, \widehat{\bld
u}_h^{j-1})
\in W_h^k\times X_h^k\times \bld V_h^k\times M_h^{k-1}$ be the numerical
solution at time $t_{j-1}$, 
and 
let $(\rho_h^{j}, \widehat{\rho}_h^{j}, \bld u_h^{j}, \widehat{\bld
u}_h^{j})
\in W_h^k\times X_h^k\times \bld V_h^k\times M_h^{k-1}$ be the numerical
solution at time $t_{j}=t_{j-1}+\delta t_{j-1}$, where $\delta t_{j-1}$ 
is the time step size at $(j-1)$-th level.
\begin{itemize}
  \item [(i)] 
Compute the maximum velocity magnitude on the mesh 
$v_{\max}:=\max_{x\in\Th} |\bld u_h^j|$, take the next time step size
\[
  \delta t_j = CFL\; h/v_{\max},
\]
and set $t_{j+1}:=t_j+\delta t_j$,
where $CFL$ is the CFL constant.
%that depends on the polynomial degree 
%$k$.
\item [(ii)] Extrapolate velocity and phase variables at time 
$t_{j+1/2}:=t_j+\frac12\delta t_j$ from data at time $t_j$ and $t_{j-1}$:
\begin{align*}
  \widetilde{\bld u}_h^{j+1/2}:=&\;\bld u_h^j
  +\frac{\delta t_j}{2\delta t_{j-1}}(\bld u_h^j-\bld u_h^{j-1}),\\
  \widetilde{\phi}_h^{j+1/2}:=&\;\phi_h^j
  +\frac{\delta t_j}{2\delta t_{j-1}}(\phi_h^j-\phi_h^{j-1}),\\
  \widetilde{\widehat{\phi}}_h^{j+1/2}:=&\;\widehat{\phi}_h^j
  +\frac{\delta t_j}{2\delta t_{j-1}}(\widehat{\phi}_h^j-\widehat{\phi}_h^{j-1}),
\end{align*}
\item [(iii)] Solve the phase field variables at the next time level 
using the extrapolated velocity:
Find $(\phi_h^{j+1}, \widehat{\phi}_h^{j+1}, \mu_h^{j+1/2},
\widehat{\mu}_h^{j+1/2})
\in 
W_h^k\times X_h^k
\times W_h^k\times X_h^k
$ such that
      \begin{subequations}
        \label{full}
        \begin{alignat}{2}
        \label{cn-phi}
        (d_t\phi_h^{j+1/2}, \psi)+
          \mathcal{C}_h^1\left((\widetilde{\bld u}_h^{j+1/2}, 
            \phi_h^{j+1/2},\widehat{\phi}_h^{j+1/2}), 
        (\psi,\widehat{\psi})\right)\quad\quad&\\
        +\mathcal{D}_h\left((\widetilde{M}_h^{j+1/2},\mu_h^{j+1/2},
          \widehat{\mu}_h^{j+1/2}), (\psi,
          \widehat{\psi})\right)
          =&\;0,\nonumber\\
        \label{full-mu}
        -(\mu_h^{j+1/2}, \eta)+
          \widetilde{\sigma}\epsilon
          (\widetilde{W}'(\phi_h^{j+1/2}), \eta) +
          \mathcal{D}_h\left((\widetilde{\sigma}\epsilon, 
            \phi_h^{j+1/2}, \widehat{\phi}_h^{j+1/2}), (\eta,
            \widehat{\eta})\right)
          =&\;0,
        \end{alignat}
      for all 
      $(\psi, \widehat{\psi},\eta,\widehat{\eta})\in 
      W^k_h\times X^k_h\times W^k_h\times X^k_h
      $, where 
      \begin{align*}
        dt\phi_h^{j+1/2}:=&\;\frac{\phi_h^{j+1}-\phi_h^{j}}{\delta
        t_j},\qquad
        \phi_h^{j+1/2}:=\;\frac12({\phi_h^{j+1}+\phi_h^{j}}),\\
        \widehat{\phi}_h^{j+1/2}:=&\;\frac12({\widehat{\phi}_h^{j+1}+
        \widehat{\phi}_h^{j}}),\quad
        \widetilde{M}_h^{j+1/2}:=\; M(\widetilde{\phi}_h^{j+1/2}),
      \end{align*}
      and 
      \[
        \widetilde{W}'(\phi_h^{j+1/2}):= 
        {W}'(\phi_h^{j})+\frac12W''(\phi_h^j)(\phi_h^{j+1}-\phi_h^j)
      \]
      is a linearization around $\phi_h^j$.
      Note that here the unknows $\phi_h^{j+1}$ and $\widehat{\phi}_h^{j+1}$
      stay at time $t_{j+1}$, while 
      the unknows $\mu_h^{j+1/2}$ and $\widehat{\mu}_h^{j+1/2}$
      stay at time $t_{j+1/2}$. The scheme \eqref{cn-phi} is linear, and can be
      efficiently implemented via static condensation such that the
      globally coupled degrees of freedom (DOFs) are those on the mesh skeleton
      only (2 DOFs per vertex, and $2(k-1)$ DOFs per dege), whose 
      computational cost is similar to a standard continuous Galerkin
      finite element method.
   %   total number 
      %and the resulting 
      %matrix sparsity pattern 
  %    is exactly the same 
  %    as the classical continuous finite elements of the same degree
  %    (after static condensation).
    \item [(iv)] Solve the velocity variables at time level $t_{j+1}$:
      Find $(\bld u_h^{j+1},\widehat{\bld u}_h^{j+1})\in \bld V_h^k\times
      M_h^{k-1}$ such that 
        \begin{align}
          \label{full-u}
          \left(\rho_h^{j+1/2}\, dt\bld u_h^{j+1/2}, \bld v\right)
          +\mathcal{C}_h^2\left((\rho_h^{j+1/2}, \widetilde{\bld
            u}_h^{j+1/2}, \widetilde{\bld u}^{j+1/2}_h),
          (\bld v,\widehat{\bld v})\right)
          \quad\qquad\quad\qquad\quad&\\
          +\mathcal{B}_h\left((\nu_h^{j+1/2}, 
            \bld u_h^{j+1/2}, 
          \widehat{\bld u}_h^{j+1/2}), (\bld v, \widehat{\bld v})
            \right)
            %\quad&\;\\
            - \mathcal{C}_h^3\left((\phi_h^{j+1/2},\widehat{\phi}_h^{j+1/2},
              \mu_h^{j+1/2}), \bld v
            \right)
          = &\;
          (\rho_h^{j+1/2}\bld f, \bld v),\nonumber
        \end{align}
      \end{subequations}
      for all 
      $(\bld v, \widehat{\bld v})\in 
      \bld V^k_h\times
      M^{k-1}_h
      $, where 
      \begin{align*}
        dt\bld u_h^{j+1/2}:=&\;\frac{\bld u_h^{j+1}-\bld u_h^{j}}{\delta
        t_j},\qquad
        \bld u_h^{j+1/2}:=\;\frac12({\bld u_h^{j+1}+\bld u_h^{j}}),\\
        \widehat{\bld u}_h^{j+1/2}:=&\;\frac12({\widehat{\bld u}_h^{j+1}+
        \widehat{\bld u}_h^{j}}),\quad
        {\rho}_h^{j+1/2}:=\; \rho({\phi}_h^{j+1/2}),\;
        {\nu}_h^{j+1/2}:=\; \nu({\phi}_h^{j+1/2}).
      \end{align*}
      In the actual implementaion, we solve for the stream function 
      $\xi_h\in \Phi_h^{k+1}$ and recover velocity via 
      the formula $\bld u_h=\nabla \times \xi_h$.
      The scheme \eqref{full-u} is can also be
      efficiently implemented via static condensation such that the
      globally coupled degrees of freedom are those on the mesh skeleton
      only (%condense out the internal DOFs of the space 
      %$\Phi_h^{k+1}$, resulting in  
      $1$ DOF per vertex, $k$ DOFs per edge
      for the space $\Phi_h^{k+1}$, and $k$ DOFs per edge for the space $M_h^{k-1}$.
    In the practical implementation, we further modify the slave variables $\rho$
    and $\nu$ as follows: 
    \[
      \rho^*(\phi):=\left\{
          \begin{tabular}{ll}
            $\rho_1$,& if $\phi>1$,
            \\[1.5ex]
            $\rho(\phi)$,& if $\phi>1$,
            \\[1.5ex]
            $\rho_2$,& if $\phi<-1$,
          \end{tabular}
        \right.
        \quad\quad\nu^*(\phi):=\left\{
          \begin{tabular}{ll}
            $\nu_1$, &if $\phi>1$,
            \\[1.5ex]
            $\nu(\phi)$,& if $\phi>1$,
            \\[1.5ex]
            $\nu_2$,& if $\phi<-1$,
          \end{tabular}
        \right.
    \]
    which ensure that 
    \[
      \min\{\rho_1,\rho_2\}\le \rho^*\le 
      \max\{\rho_1, \rho_2\}, \text{ }
      \min\{\nu_1,\nu_2\}\le \nu^*\le 
      \max\{\nu_1, \nu_2\}.
    \]
\end{itemize}

      \section{Numerical results}
\label{sec:num}
In this section, we present several numerical results for the 
divergence-free HDG scheme \eqref{full}
proposed in Section \ref{sec:scheme}.
The NGSolve software \cite{Schoberl16} is used for the simulation.

\subsection{Accuracy test}
We use the method of manufactory solutions to test the spatial and temporal
accuracy of the scheme. In particular, we take the computational domain to
be a periodic unit square, add a source term $f_\phi$ in the equation 
\eqref{ns-eq1}, and take the source terms $f_\phi$ in \eqref{ns-eq1} and 
$\bld f$ in \eqref{ns-eq2} such
that the exact solution is given  as follows:
\begin{align*}
  \phi(t,x) =&\; \sin(\pi t)\sin(2\pi x)\sin(2\pi y),\\
  \mu(t,x) =&\; \widetilde{\sigma}(\epsilon^{-1}W'(\phi(t,x))-\epsilon\triangle
  \phi(t,x)),\\
  \bld u(t,x) =&\;\left[
  \begin{tabular}{l}
    $0.2\sin(\pi t)\sin(2 \pi x)\cos(2\pi y)$\\
    $-0.2\sin(\pi t)\cos(2 \pi x)\sin(2\pi y)$
\end{tabular}\right].
\end{align*}
Furthermore, we take $\rho_1=100$, $\rho_2=10$, $\nu_1=10$, $\nu=1$, 
$\sigma=10$, $\epsilon=0.04$, and $\gamma= 10^{-3}\epsilon$.
The final time is taken to be $T=0.5$, and we use a uniform time step size 
$\delta t = h^{(k+1)/2}$. The history of convergence of the $L^2$-norm
errors in $\phi_h$, $\mu_h$, and $\bld u_h$ at the final time $T=0.5$
on a sequence of
uniform square meshes is recorded in Table \ref{table:err} for $k=1$ and 
$k=2$. We clearly observe optimal convergence order of $k+1$ 
in all the variables for both cases, except for the error in $\mu_h$ when
$k=1$, where we loose half order of convergence. 
In particular, this indicates the expected second-order accuracy in time of the
%CNAB
temporal discretization.
%Hence, the scheme \eqref{full} is observed to be 
%second-order accurate in time, and $(k+1)$-th order accurate in space.

\begin{table}[ht!]
\begin{center}
% \footnotesize
%\scalebox{0.78}
{%
\begin{tabular}{| c|c | c | c | c | c | c | c |
} \hline
  &   &
  \multicolumn{2}{c|}{$\|u-u_h\|$} & \multicolumn{2}{c|}{$\|\phi-\phi_h\|$}& 
  \multicolumn{2}{c|}{$\|\mu-\mu_h\|$} 
  \\\hline
  $k$& $1/h$ & Error & Order  & Error & Order  & Error & Order
  %& Error & Order
     \\
\hline
 & 8 & 5.68e-03  &  -- & 1.96e-02  &  -- & 1.13e+01  &  -- \\                    
& 16 & 1.10e-03  &  2.37 & 5.80e-03  &  1.75 & 2.62e+00  &  2.11 \\              
1& 32 & 2.62e-04  &  2.07 & 1.54e-03  &  1.91 & 7.29e-01  &  1.84 \\              
& 64 & 6.41e-05  &  2.03 & 3.92e-04  &  1.98 & 2.27e-01  &  1.68 \\             
& 128 & 1.59e-05  &  2.01 & 9.84e-05  &  1.99 & 7.88e-02  &  1.53 \\    
\hline
& 8 & 4.39e-04  &  -- & 1.38e-03  &  -- & 1.50e+00  &  -- \\ 
& 16 & 5.08e-05  &  3.11 & 2.05e-04  &  2.75 & 2.26e-01  &  2.73 \\
2& 32 & 6.05e-06  &  3.07 & 2.47e-05  &  3.06 & 2.83e-02  &  2.99 \\
& 64 & 7.58e-07  &  3.00 & 3.16e-06  &  2.96 & 3.52e-03  &  3.01 \\ 
& 128 & 9.48e-08  &  3.00 & 4.04e-07  &  2.97 & 4.38e-04  &  3.01 \\
\hline
\end{tabular}
  }
\caption{\it History of convergence of the $L^2$-errors.}
\label{table:err}
\end{center}
\end{table}

\subsection{Rising bubble}
\begin{figure}[h!]
\centering
\includegraphics[width=.25\textwidth]{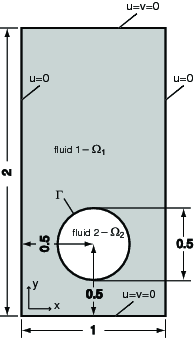}
  \caption{Initial configuration and boundary conditions for the rising bubble
  problem.}
\label{fig:bubble}
\end{figure}

We consider the rising bubble benchmark problem proposed in
\cite{Hysing09}. 
The test setup is extensively described in \cite{Hysing09}. The domain 
$\Omega=[0,1]\times[0,2]$ is filled with fluid 1 $(\phi\approx 1)$ except
for a circular bubble, which consists of fluid 2 $(\phi\approx -1)$. The
initial bubble has a radius of $0.25$ with its center at $(0.5,0.5)$.
See Figure~\ref{fig:bubble} for the sketch of the domain and boundary conditions.

Table \ref{table:bubble} 
lists the fluid and physical parameters for the two test cases used in
\cite{Hysing09}. 
%\begin{wraptable}{l}{.5\textwidth}
\begin{table}[ht!]
\begin{center}
  %\footnotesize
  %\scalebox{1}
  {
  \begin{tabular}{ c  c  c  c c c c c} \hline
    Test Case & $\rho_1$ &$\rho_2$& $\mu_1$&$\mu_2$&$\bld f$
              &$\sigma$ \\
    1 &	1000 &	100 	&10 &1 &	(0,-0.98)& 	24.5\\
    2 &	1000 &	1 	&10 	&0.1 &	(0,-0.98)& 	1.96\\
    \hline 
\end{tabular}
}
\end{center}
\caption{\it Physical parameters 
  %and final time 
for the test cases.}
\label{table:bubble}
%\end{wraptable}
\end{table}

Because the density of the bubble is smaller than the density of the
surrounding fluid ($\rho_2< \rho_1$),
the bubble rises. The evolution of the bubble is tracked for three time units during which the defined
benchmark quantities are measured.
The measured benchmark quantities are the following:
\begin{itemize}
\item Center of mass
  \[
    y_c = \frac{\int_{\phi<0} y\,\mathrm{dx}}{
    \int_{\phi<0} 1\,\mathrm{dx}},
  \]
  where $y$ is the vertical coordinate.
\item Circularity:
  \[
    c=\frac{\text{perimeter of area-equivalent circle}}{\text{perimeter of
    bubble}}=\frac{2\sqrt{\int_{\phi<0}\pi\,\mathrm{dx}}}{\int_{\phi=0}1\,\mathrm{ds}}.
  \]
\item Rise velocity: 
  \[
    V_c = \frac{\int_{\phi<0} v\,\mathrm{dx}}{
    \int_{\phi<0} 1\,\mathrm{dx}}
  \]
  where $v$ is the vertical component of the velocity $\bld u$.
\end{itemize}

For both test cases,
the viscous effects and surface tension forces dominates the convection.
We take polynomial degree $k=2$ and use a uniform time step size throughout. 
Uniform square mesh with mesh size $h=2^{-5}$, $h=2^{-6}$, and $h=2^{-7}$ 
are used in the simulation. We take time step size to be $\delta t = 0.005$
(600 total steps) for  the cases $h=2^{-5}$ and $h=2^{-6}$, and $\delta t =
0.0025$ (1200 total steps) for the case $h=2^{-7}$.
Following \cite{AlanVoigt12}, we further take 
the diffuse interface width  $\epsilon=0.64 h$, and mobility coefficient 
$\gamma=10^{-3}\epsilon$.
\subsubsection{Results for test case 1}
\begin{figure}[h!]
\centering
\includegraphics[width=.45\textwidth]{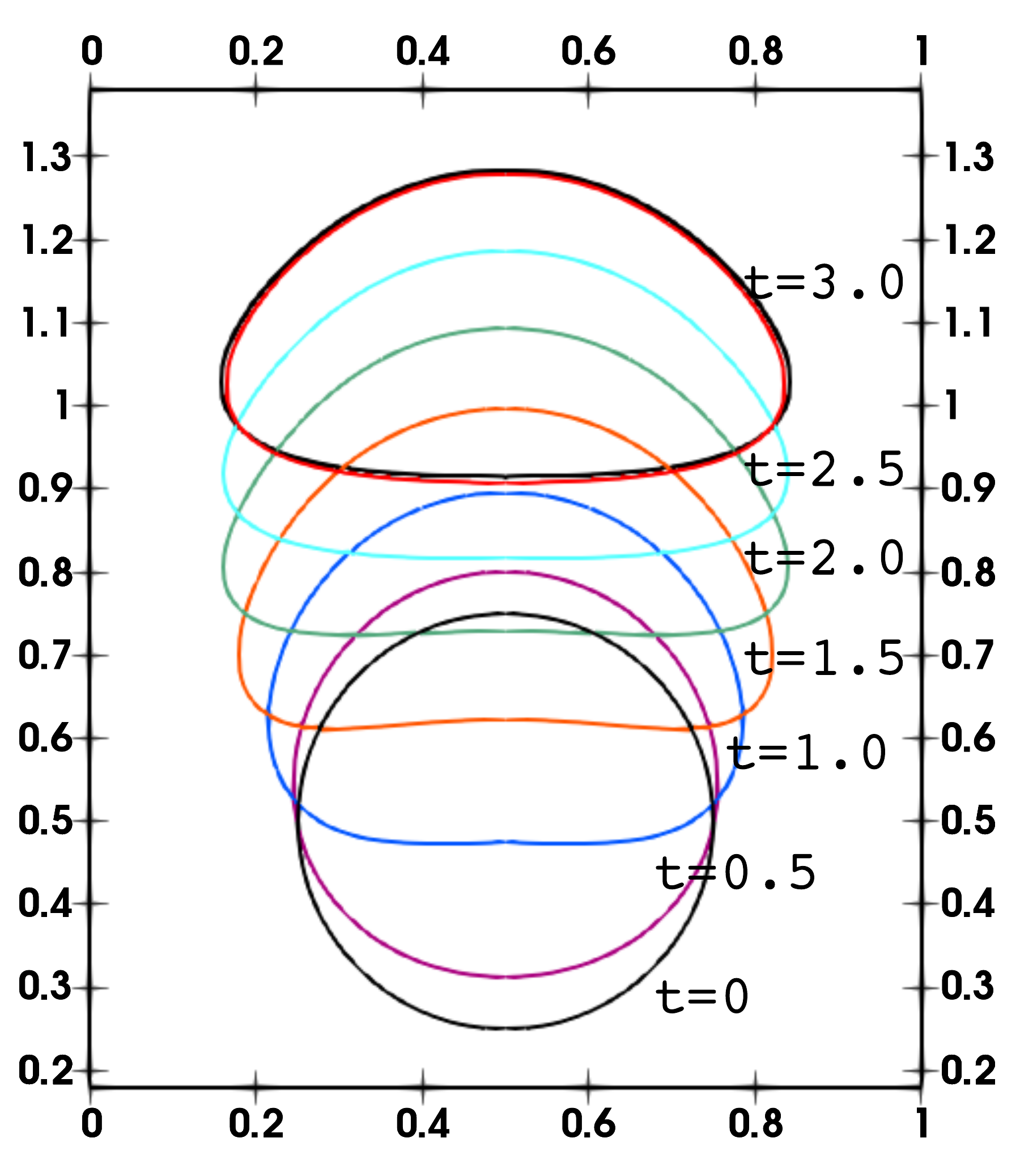}
\caption{Shapes of the rising bubble 
(contour line of $\phi_h=0$) 
  at  various time
  %time  $t=0,0.5,1,1.5,2,2.5,3$. 
  for $h=2^{-6}$ and the reference shape at $t=3$ (in black).
 % Polynomial degree $k=2$, time step size $\delta t = 0.005$.
(For interpretation of the colors in this figure, the reader is referred to the web version of this article.)}
\label{fig:bs2}
\end{figure}

The bubble, being initially circular, first stretches horizontally and 
develops a dimple at the bottom
before it reaches a stable ellipsoidal shape, see Figure~\ref{fig:bs2}.

\begin{figure}[h!]
\centering
\includegraphics[width=.6\textwidth]{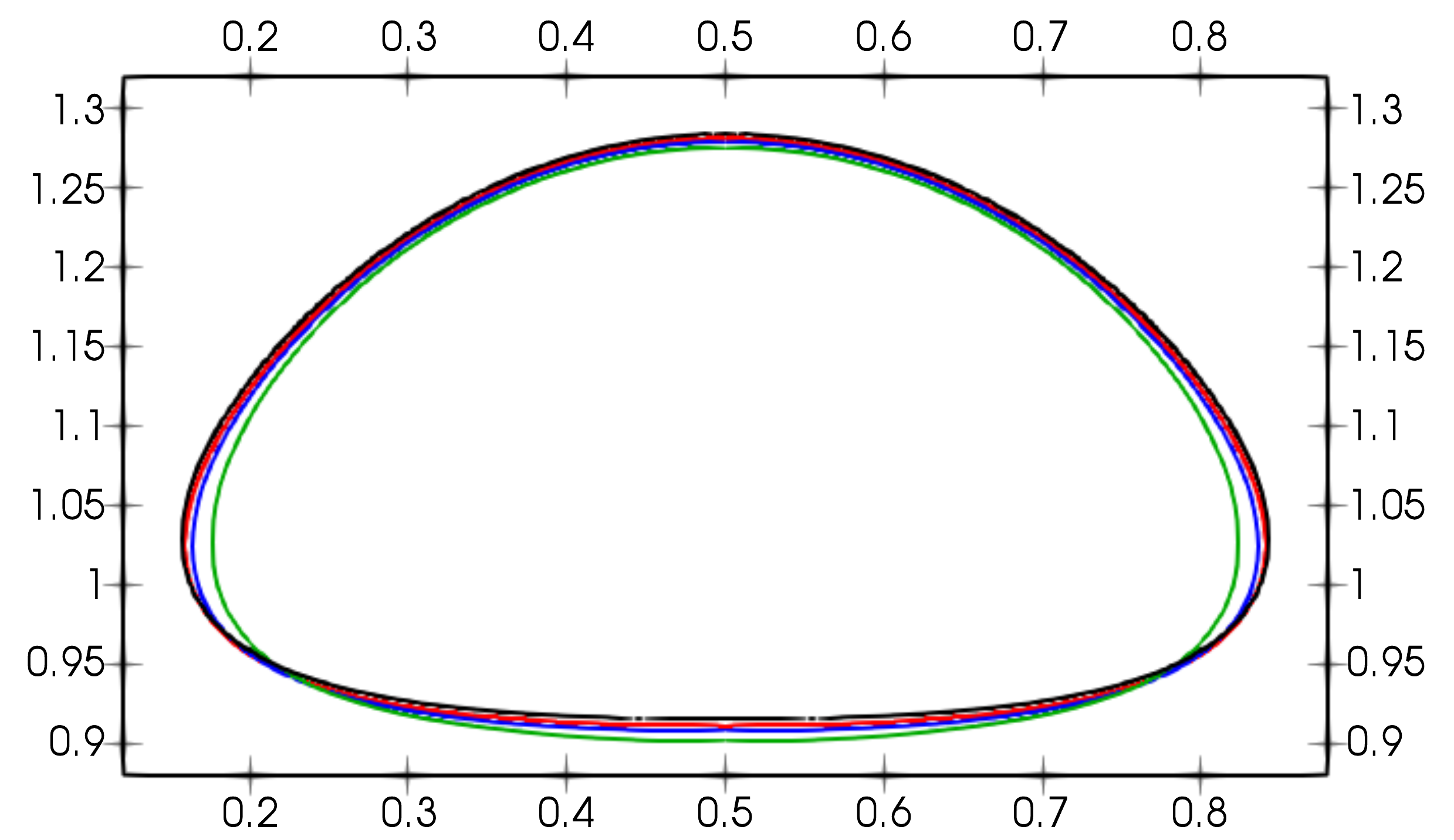}
  \caption{Shape of the rising bubble at final time $t=3$ 
    for different $h$. 
    %Polynomial degree $k=2$. 
  Green: $h=2^{-5}$. Blue: $h=2^{-6}$. Red: $h=2^{-7}$. Black:
  reference data from \cite{Hysing09}.
(For interpretation of the colors in this figure, the reader is referred to the web version of this article.)}
\label{fig:bs}
\end{figure}

Figure \ref{fig:bs} shows the bubble shapes 
(contour line of $\phi_h=0$) 
at the final time ($t=3$),
along with reference data from \cite{Hysing09}. It is clear that as mesh
size $h$ (and $\epsilon$) decreases, the bubble shape converges to the reference value. 

Table~\ref{table:bb} show the quantitative comparison with the benchmark
values. The line 'ref' gives a reference value from group $3$ in
\cite{Hysing09}. One can see
that each quantity approaches the reference value as  $h$ decreases.
\begin{table}[ht!]
\begin{center}
  %\footnotesize
  \scalebox{1}
  {
  \begin{tabular}{ c  c  c  c c c c c} 
    \hline
    $h$ & $c_{\min}$ &$t|_{c=c_{\min}}$& $V_{c,\max}$&$t|_{V_c=V_{c,\max}}$&
    $y_c(t=3)$ \\
    \hline
    $2^{-5}$ &  0.9166 &   1.885  &  0.2371 &   0.970 &   1.0732\\
    $2^{-6}$ &  0.9067 &   1.905  &  0.2388 &   0.930 &   1.0754\\
    $2^{-7}$ &  0.9034 &   1.895  &  0.2400 &   0.925 &   1.0776\\
    ref      &  0.9013 &   1.900  &  0.2417 &   0.924 &   1.0799\\
    \hline 
\end{tabular}
}
\end{center}
\caption{\it 
 Minimum circularity and maximum rise velocity, with corresponding
incidence times and final position of the center of mass for test case 1.}
\label{table:bb}
%\end{wraptable}
\end{table}

\begin{figure}[h!]
\centering
\includegraphics[width=.32\textwidth]{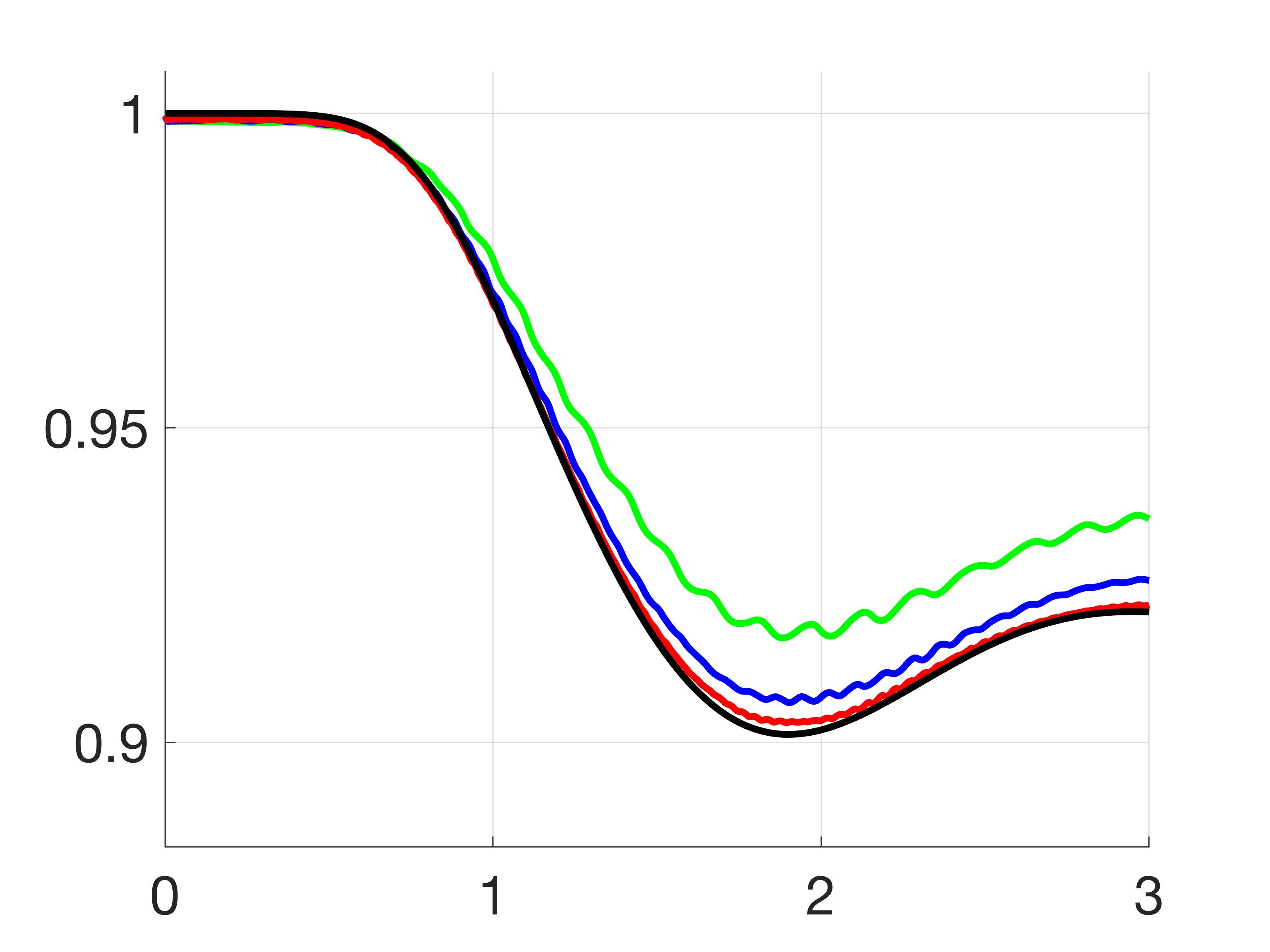}
\includegraphics[width=.32\textwidth]{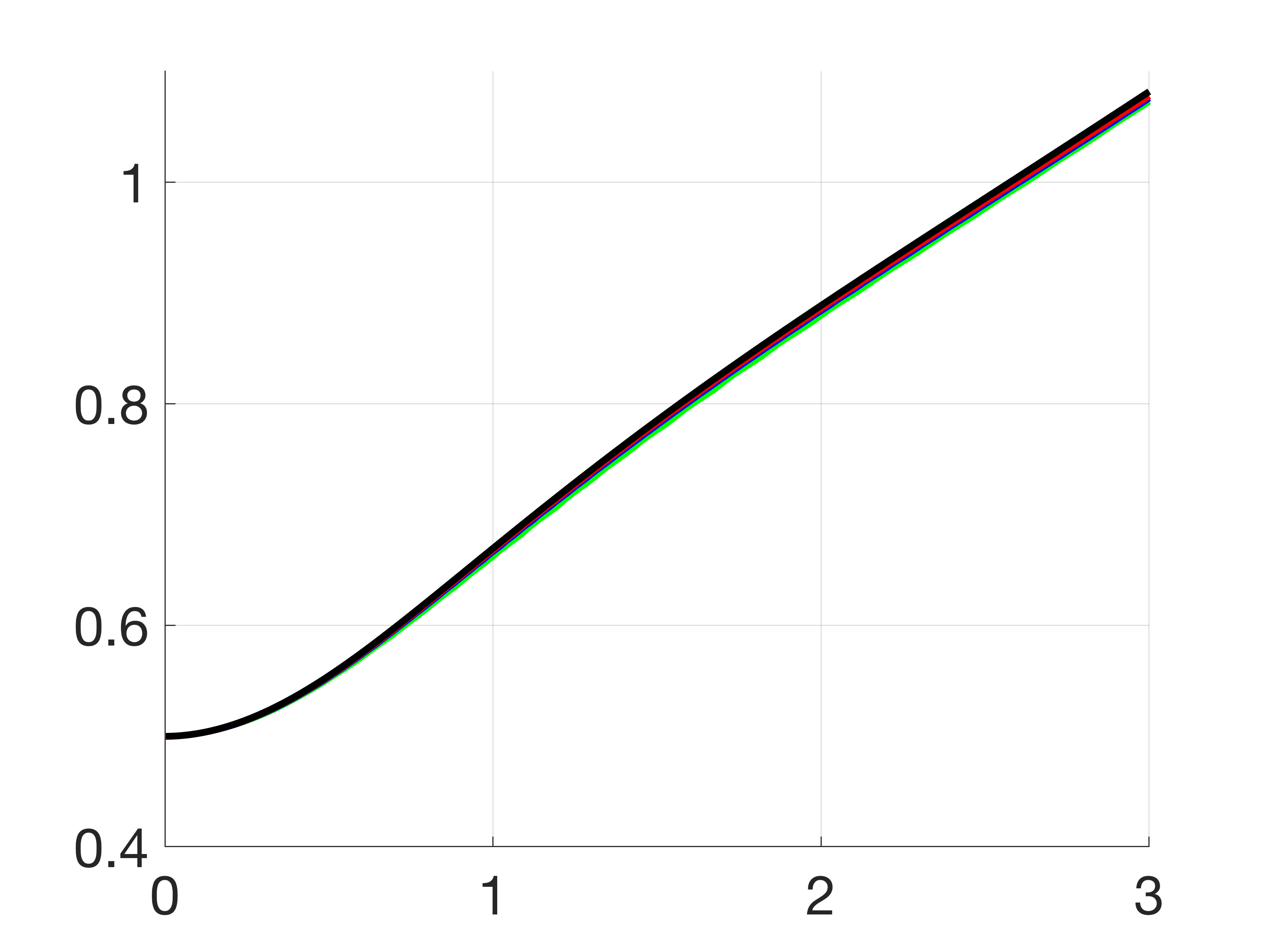}
\includegraphics[width=.32\textwidth]{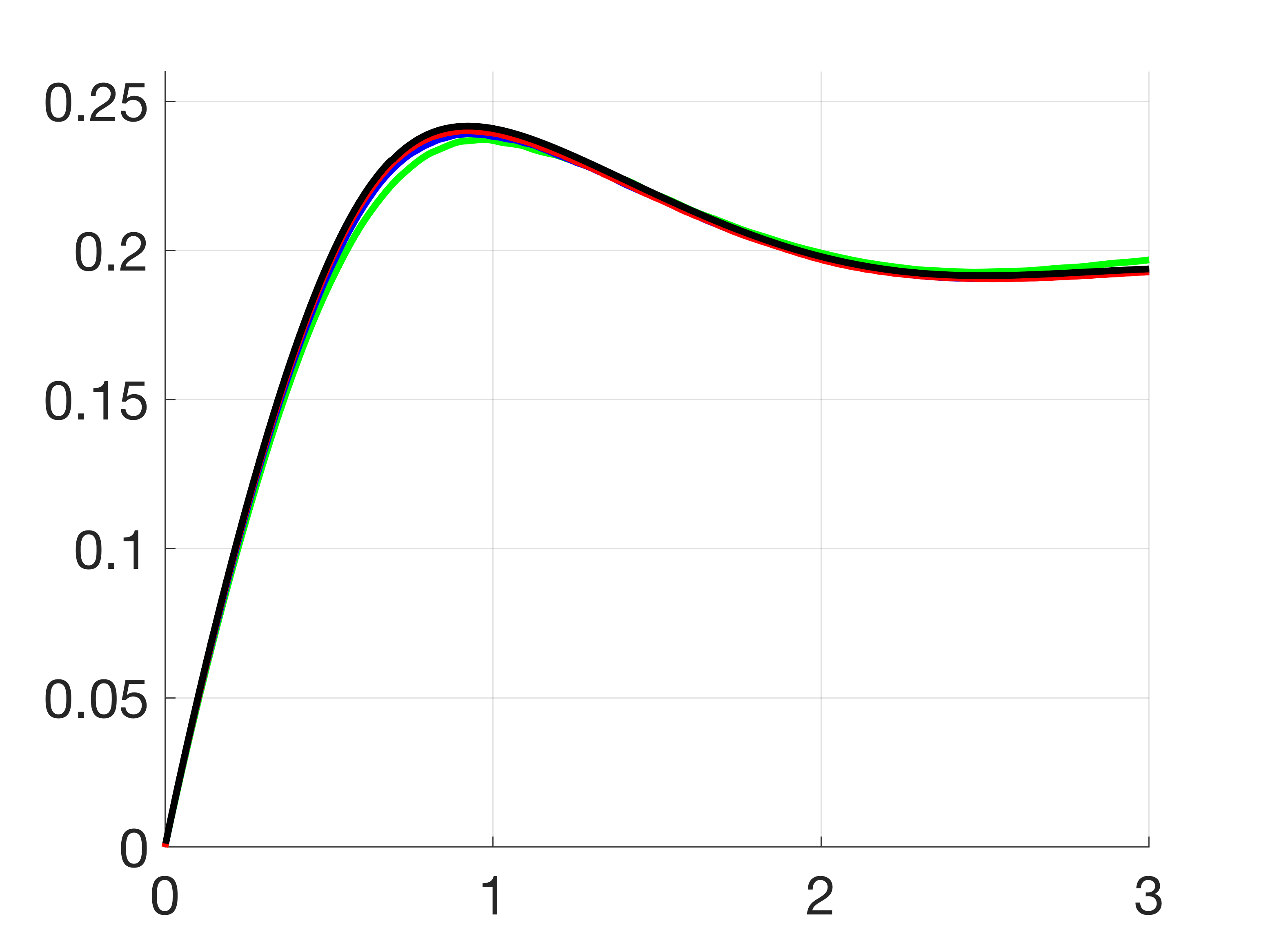}
  \caption{
    Circularity (left), center of mass (middle), and rising velocity
    (right) over time for test case 1.
  Green: $h=2^{-5}$. Blue: $h=2^{-6}$. Red: $h=2^{-7}$. 
%  Polynomial degree $k=2$.
  Black:
  reference data from group 3 of \cite{Hysing09}.
(For interpretation of the colors in this figure, the reader is referred to the web version of this article.)}
\label{fig:bb3}
\end{figure}

Furthermore, we plot the circularity, center of mass, and rise velocity over time 
in Figure~\ref{fig:bb3}.
All the quantities seem to converge as mesh size $h$ decreases.

\subsubsection{Results for test case 2}
\begin{figure}[h!]
\centering
\includegraphics[width=.23\textwidth]{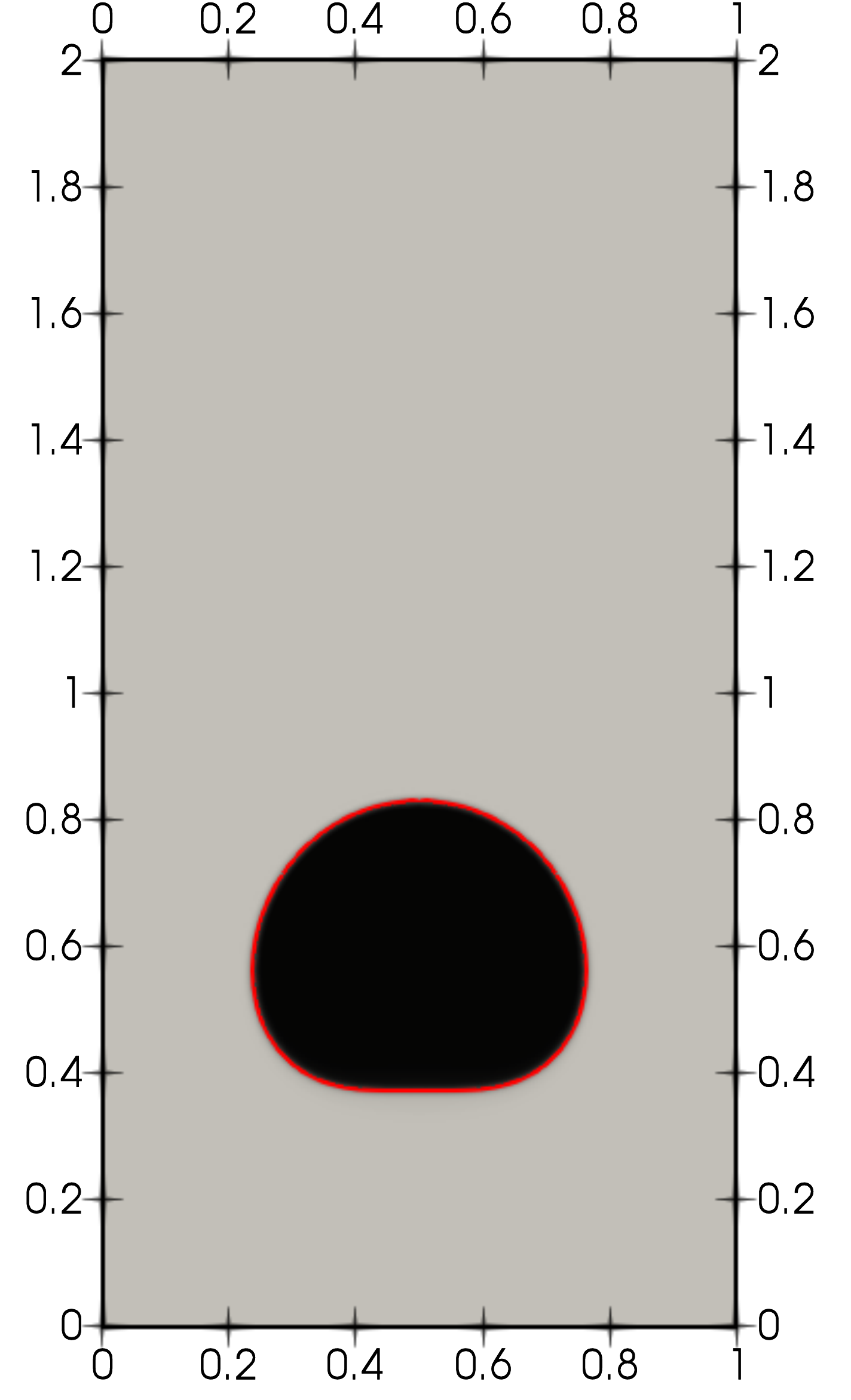}
\includegraphics[width=.23\textwidth]{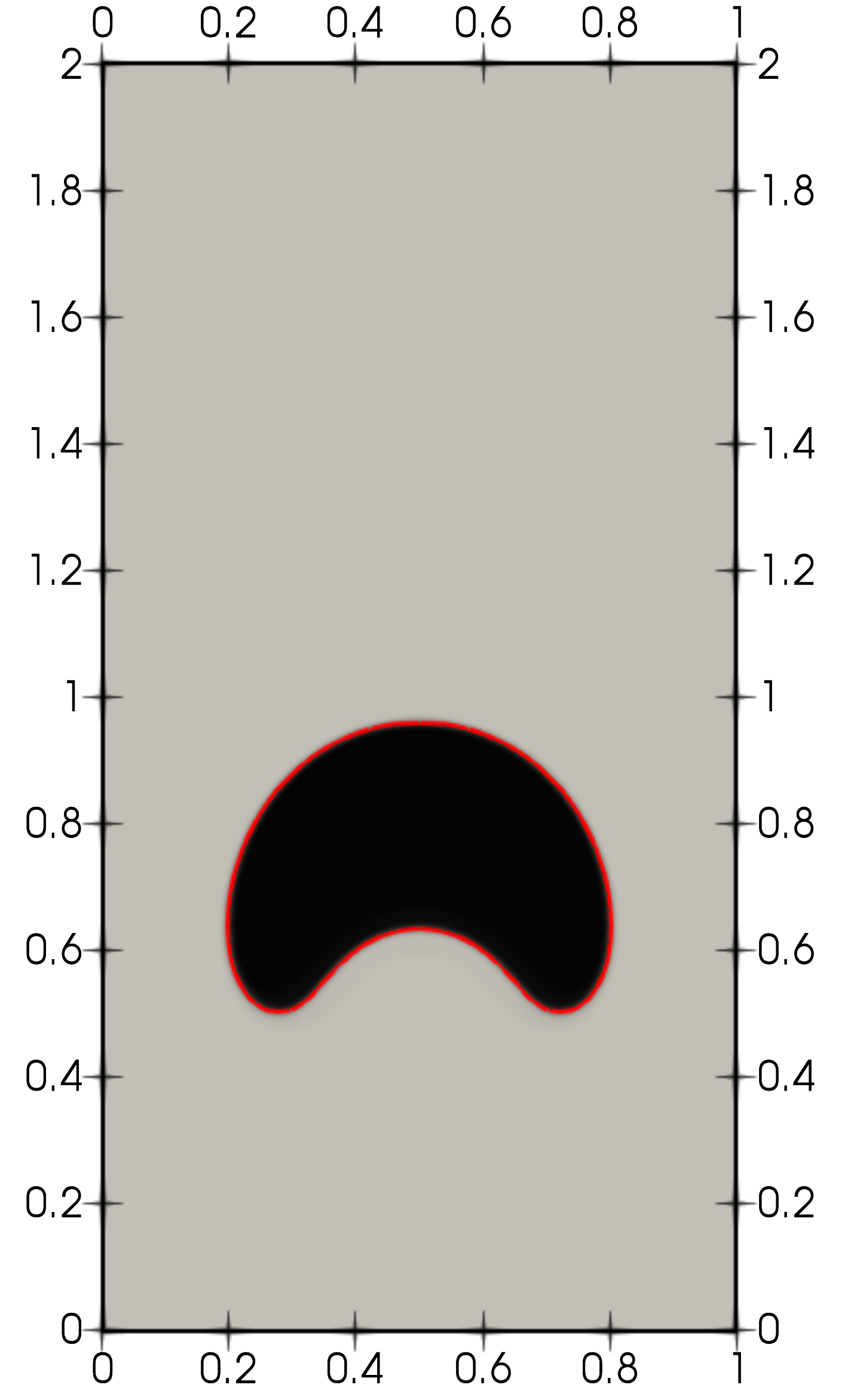}
\includegraphics[width=.23\textwidth]{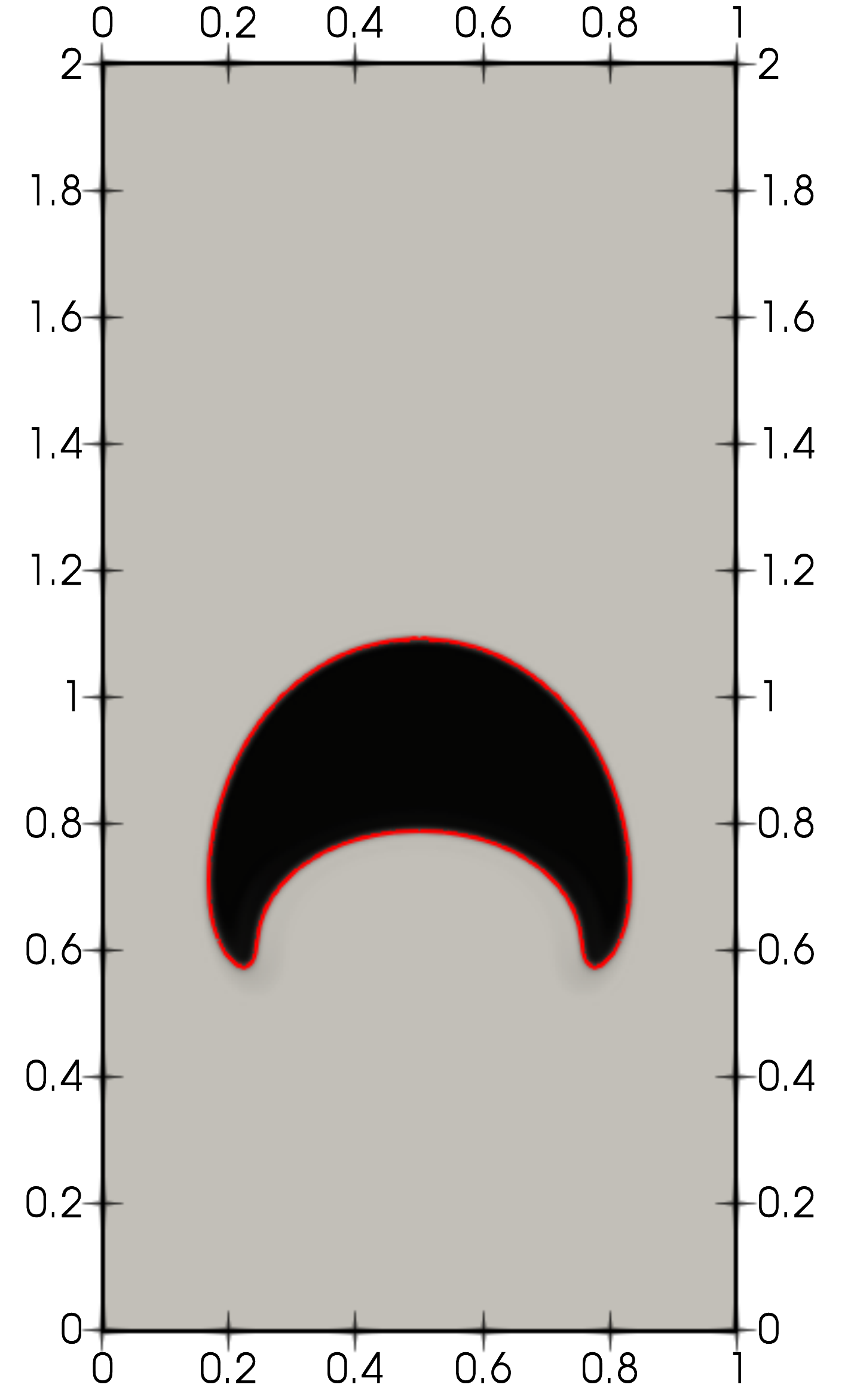}
\includegraphics[width=.23\textwidth]{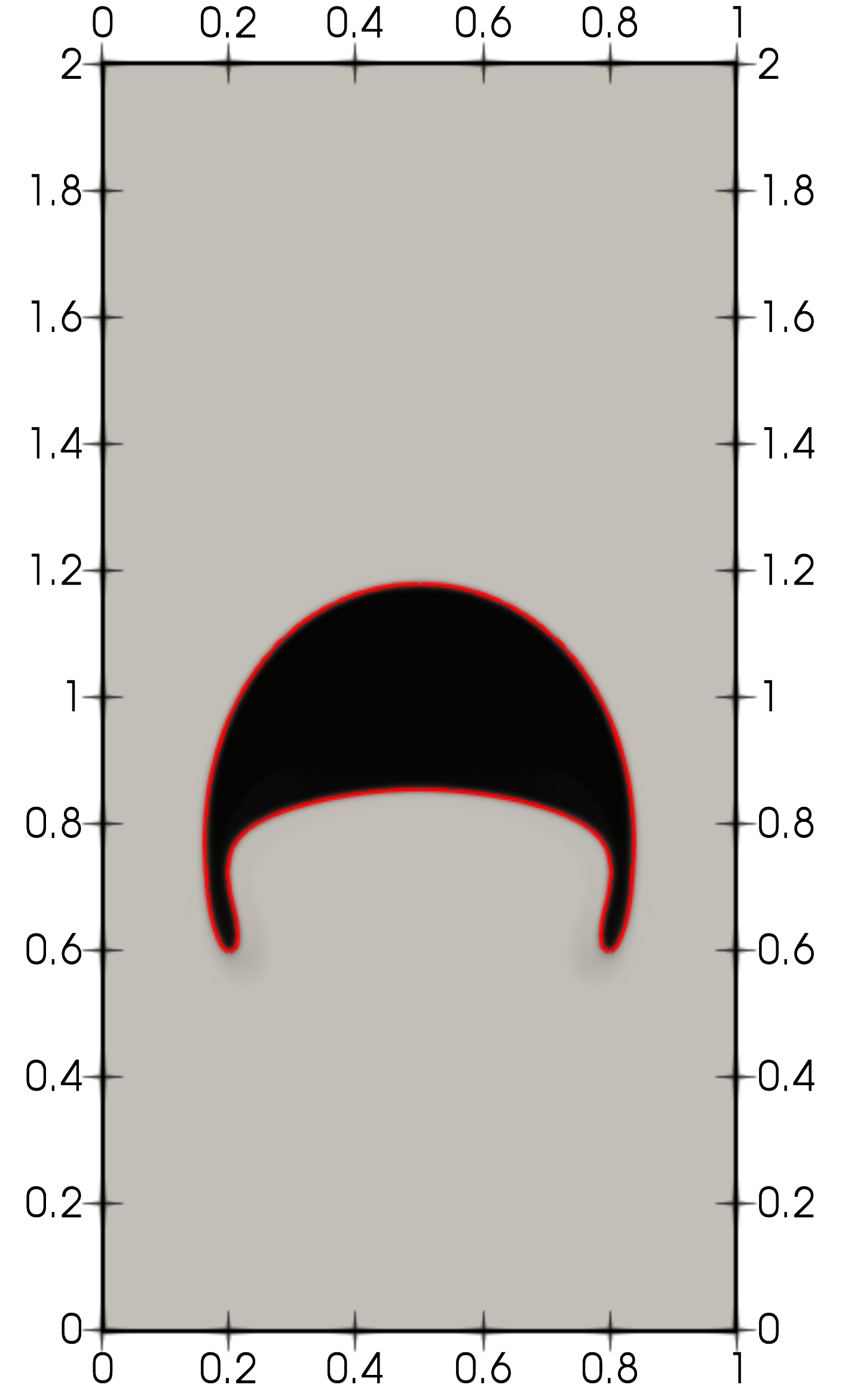}
\includegraphics[width=.23\textwidth]{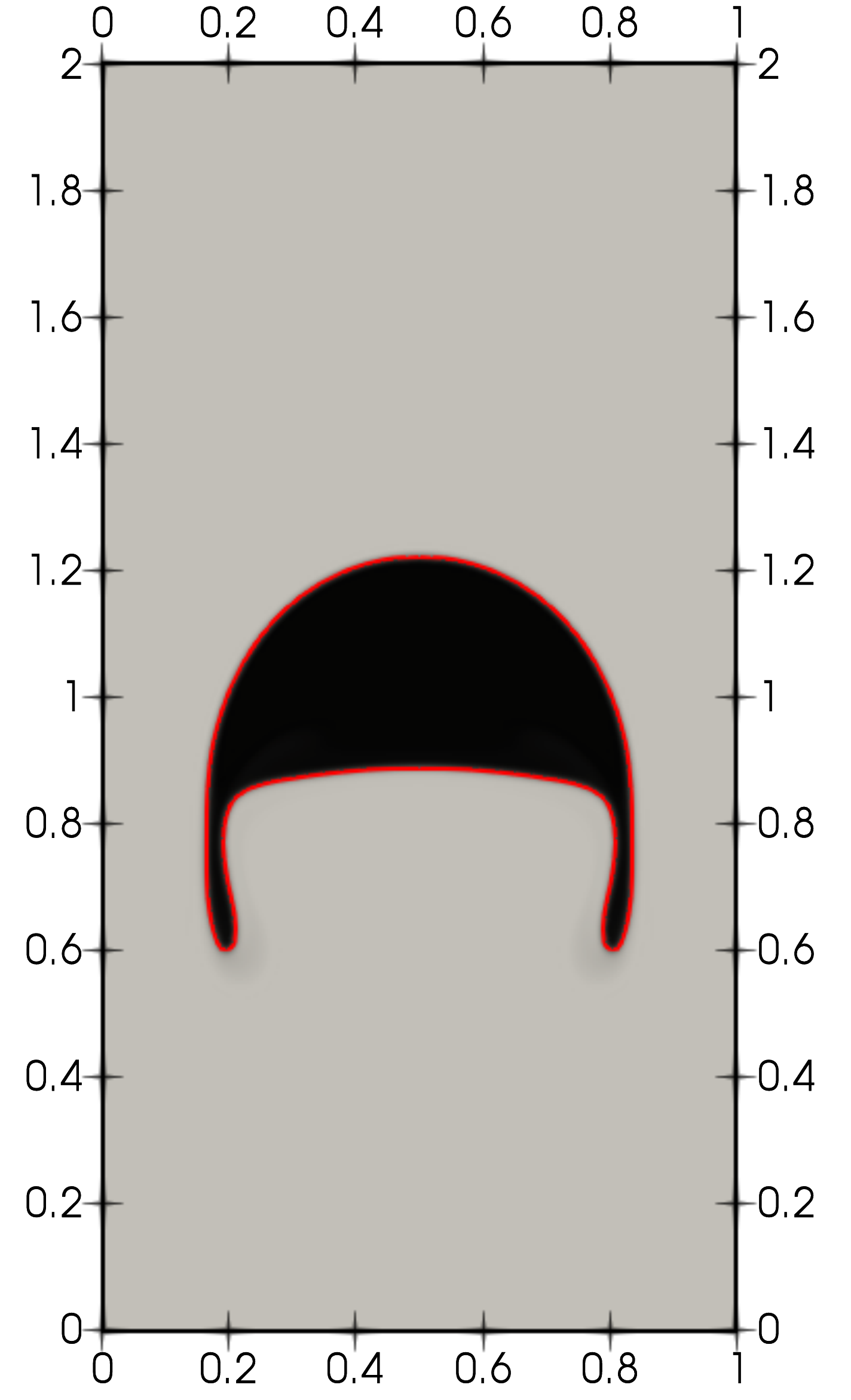}
\includegraphics[width=.23\textwidth]{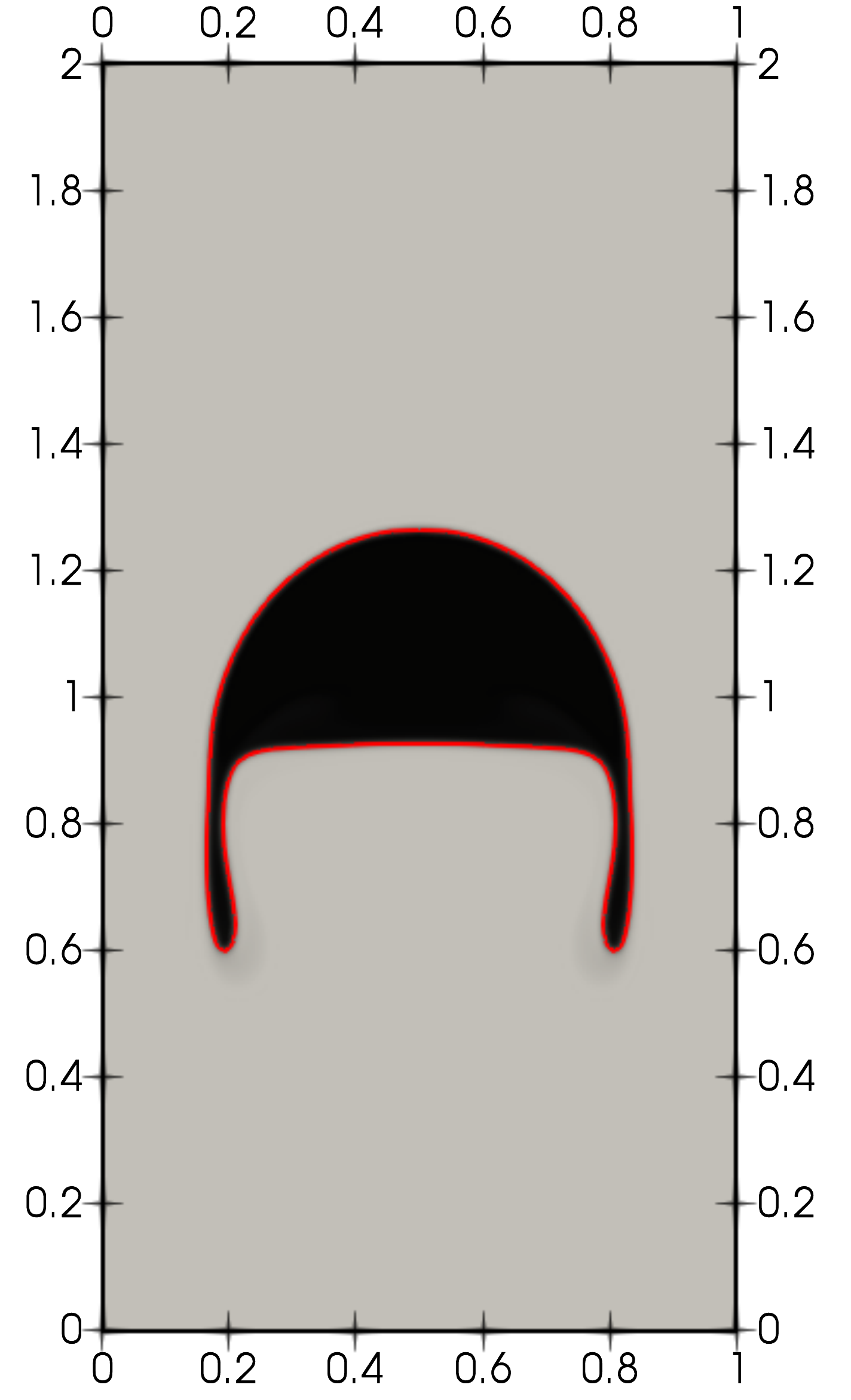}
\includegraphics[width=.23\textwidth]{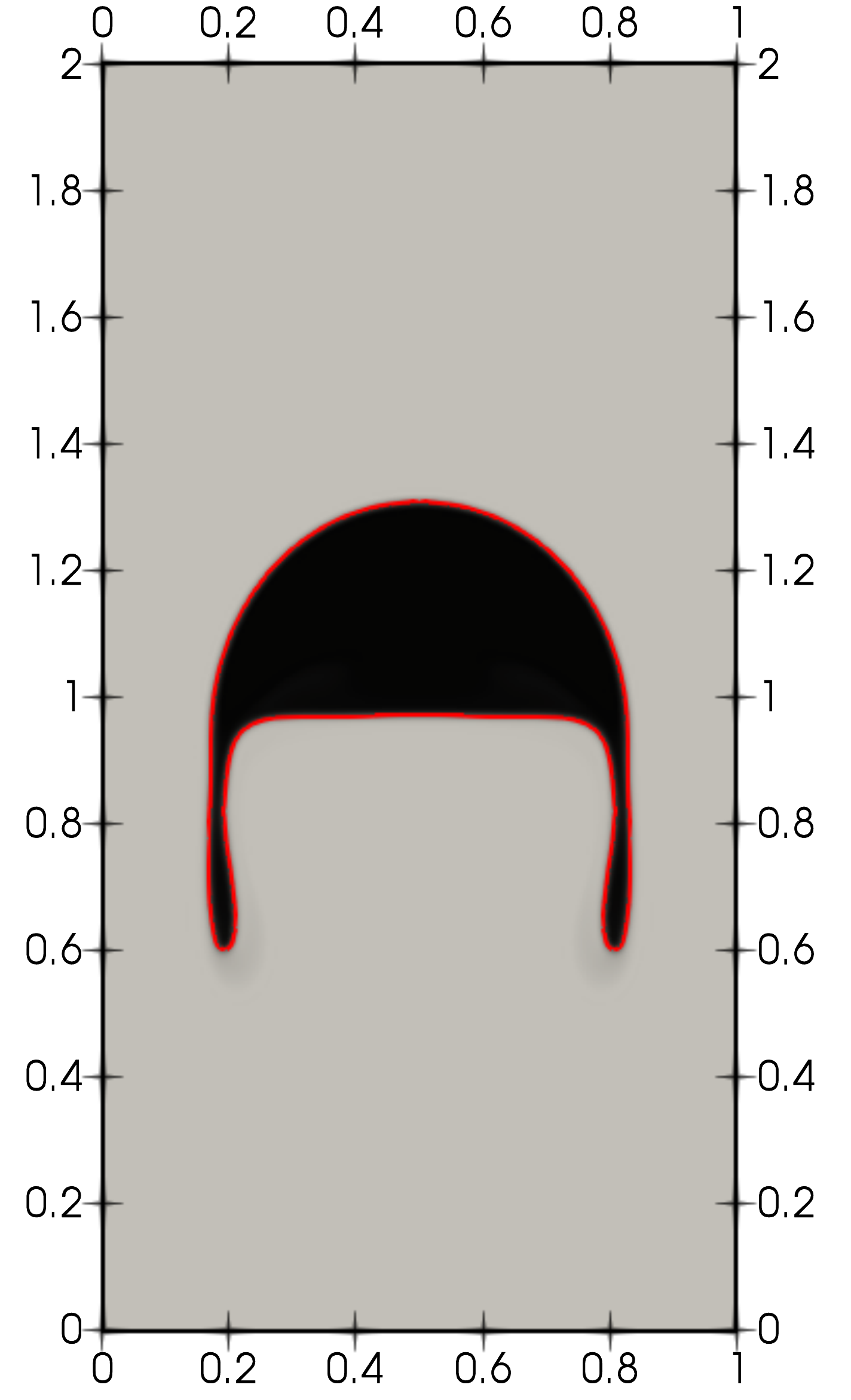}
\includegraphics[width=.23\textwidth]{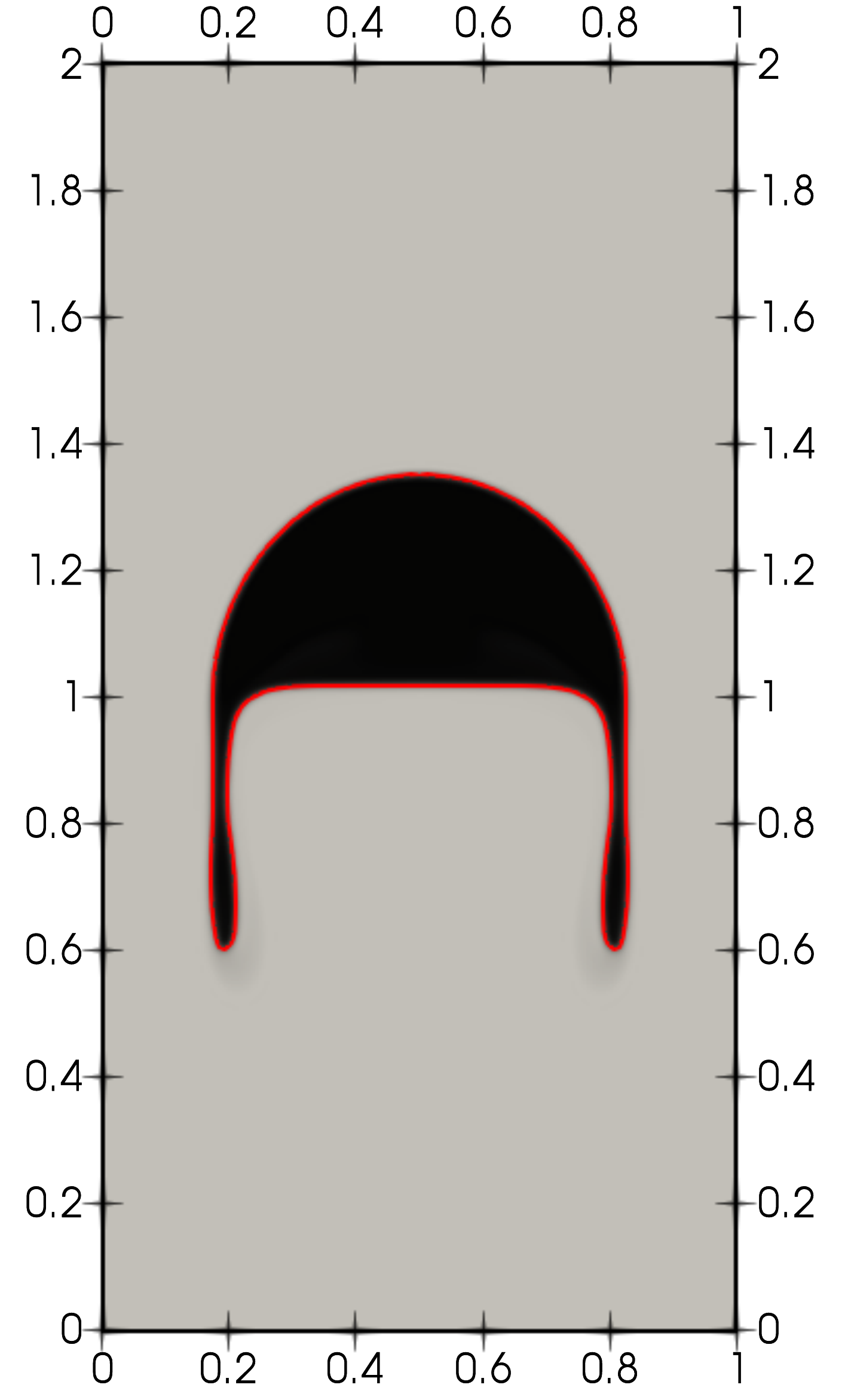}
\caption{Contour of the rising bubble for test case 2 
at time  $t=0.6,1.2,1.8,2.2,2.4,2.6,2.8, 3.0$ (from left to right, top to
bottom)
  for $h=2^{-7}$. 
  Red contour: the bubble interface $\phi_h=0$.
(For interpretation of the colors in this figure, the reader is referred to the web version of this article.)}
\label{fig:bbt2}
\end{figure}

In test case 2, the decrease in surface tension causes the bubble to develop a more non-convex shape
and thin filaments.
%, which eventually break off. 
This is a much harder problem, and also in \cite{Hysing09},
agreement between the used numerical approaches could not be achieved. It even remains unclear
if break off of the thin filaments should occur for this setting. 
%The results in [1] show that the size of the filaments
%might become very small (e.g., 0.01). Therefore, one can expect phase field methods to have some
%trouble because this can be close to the interface thickness 	.
We refer to Figure~\ref{fig:bbt2} for a depiction of the 
temporal evolution of the bubble shape at mesh size $h=2^{-7}$.
The initial circular shape is shown to gradually develop two filaments on its sides 
as it experiences an upward pushing force. But the filaments does not
break.
This observation is similar to that for the diffusive interface method used in 
\cite{AlanVoigt12}.

\begin{figure}[h!]
\centering
\includegraphics[width=.5\textwidth]{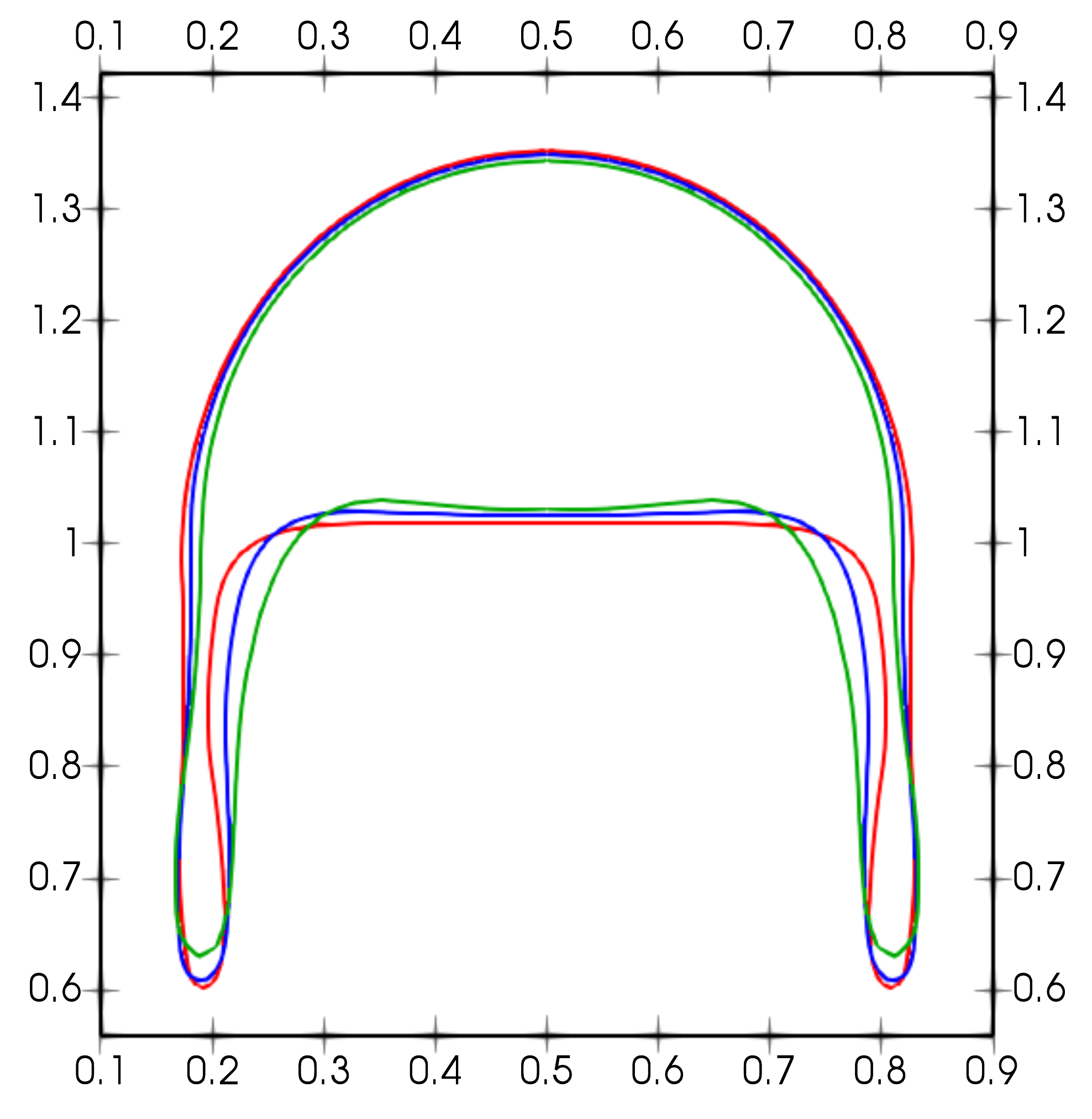}
  \caption{Shape of the rising bubble for test case 2 at final time $t=3$ 
  for different $h$. 
  Green: $h=2^{-5}$. Blue: $h=2^{-6}$. Red: $h=2^{-7}$. 
(For interpretation of the colors in this figure, the reader is referred to the web version of this article.)}
\label{fig:bbt2x}
\end{figure}

Figure~\ref{fig:bbt2x} 
show the bubble shapes 
at the final time ($t=3$).  
One can see that, the
filaments become thinner for smaller  $h$ (and smaller $\epsilon$). 
Therefore, it is possible that break off happens in the case 
$\epsilon\rightarrow
0$. 

Similar to test case 1, we list in Table~\ref{table:bbt2} 
the quantitative comparison with the benchmark
values. 
Here, we restricted the comparison with the time interval $[0,2]$ because also the reference solutions
do not agree well for later times.
One can again see
that each quantity approaches the reference value as  $h$ decreases.
\begin{table}[ht!]
\begin{center}
  %\footnotesize
  \scalebox{1}
  {
  \begin{tabular}{ c  c  c  c c c c c} 
    \hline
    $h$ & $c_{\min}$ &$t|_{c=c_{\min}}$& $V_{c,\max}$&$t|_{V_c=V_{c,\max}}$&
    $y_c(t=2)$ \\
    \hline
    $2^{-5}$ &  0.6629  &  2.0000 &   0.2491 &   0.7500  &  0.9026\\
    $2^{-6}$ &0.6627  &  2.0000 &   0.2489 &   0.7400  &  0.9070\\
    $2^{-7}$ &0.6678  &  2.0000 &   0.2494 &   0.7275  &  0.9107 \\
    ref      &0.6901 &2.0000 &0.2502 &0.7300 &0.9154\\
    \hline 
\end{tabular}
}
\end{center}
\caption{\it 
 Minimum circularity and maximum rise velocity, with corresponding
incidence times and final position of the center of mass for test case 2.}
\label{table:bbt2}
%\end{wraptable}
\end{table}

\begin{figure}[h!]
\centering
\includegraphics[width=.32\textwidth]{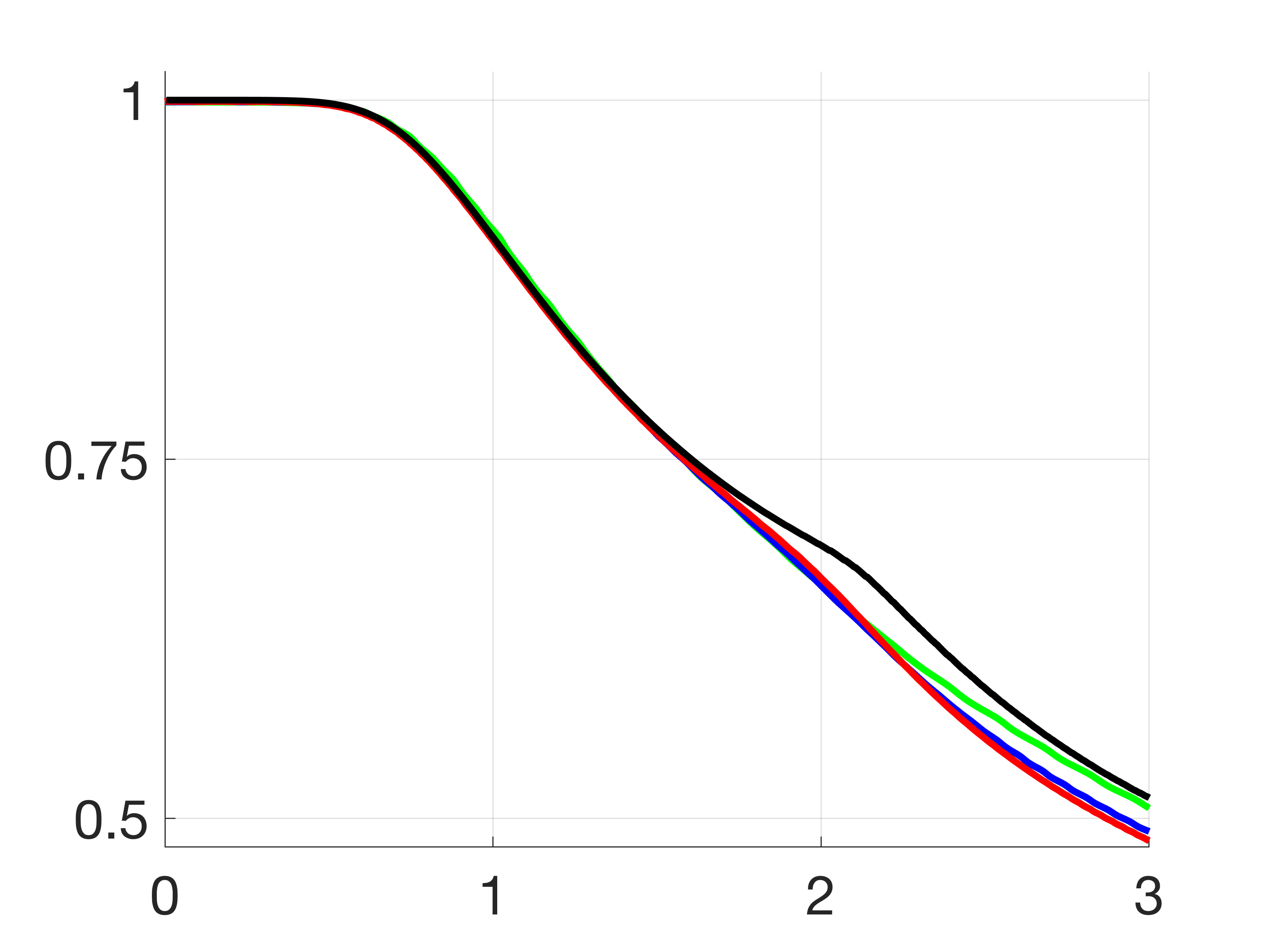}
\includegraphics[width=.32\textwidth]{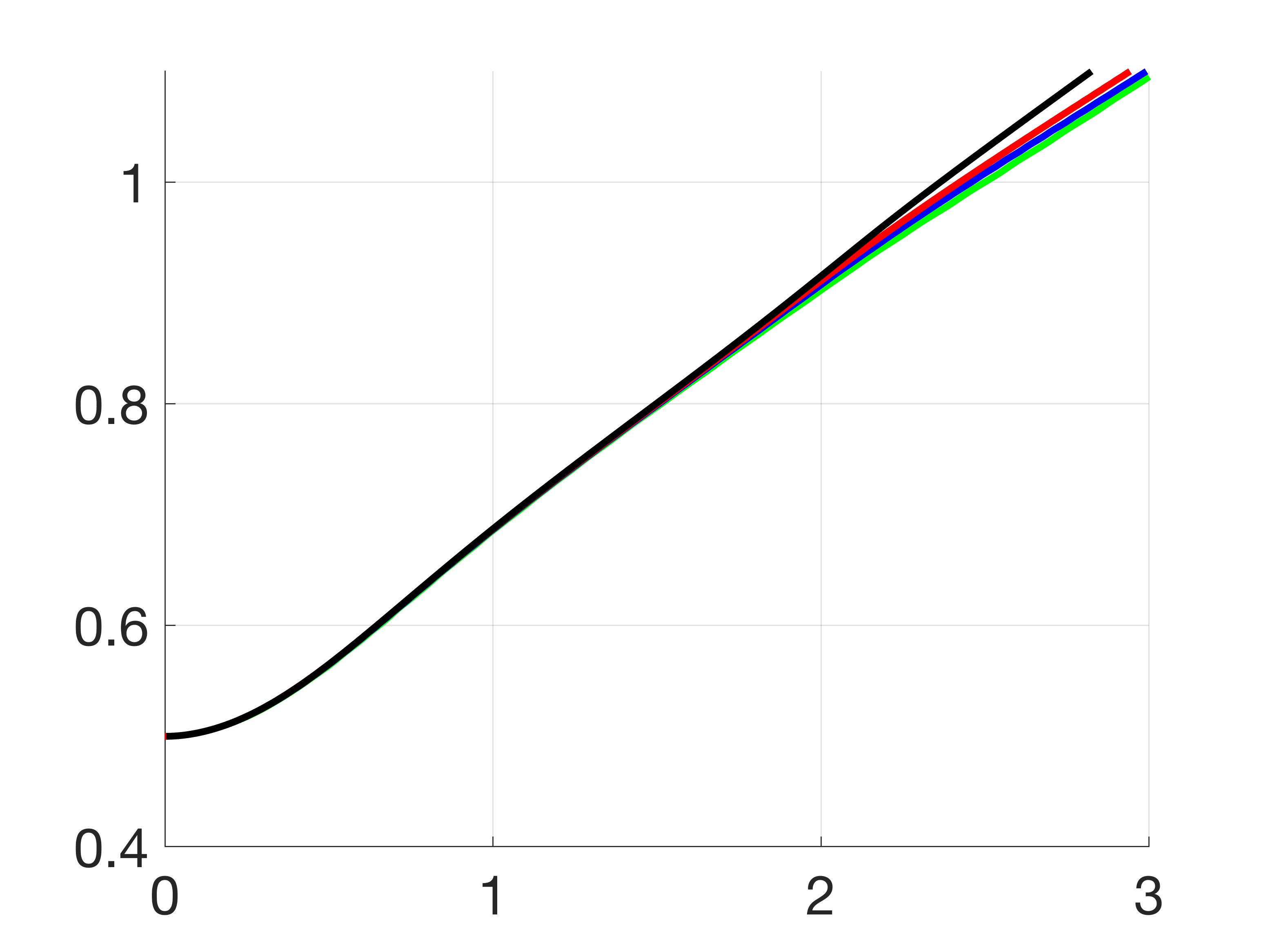}
\includegraphics[width=.32\textwidth]{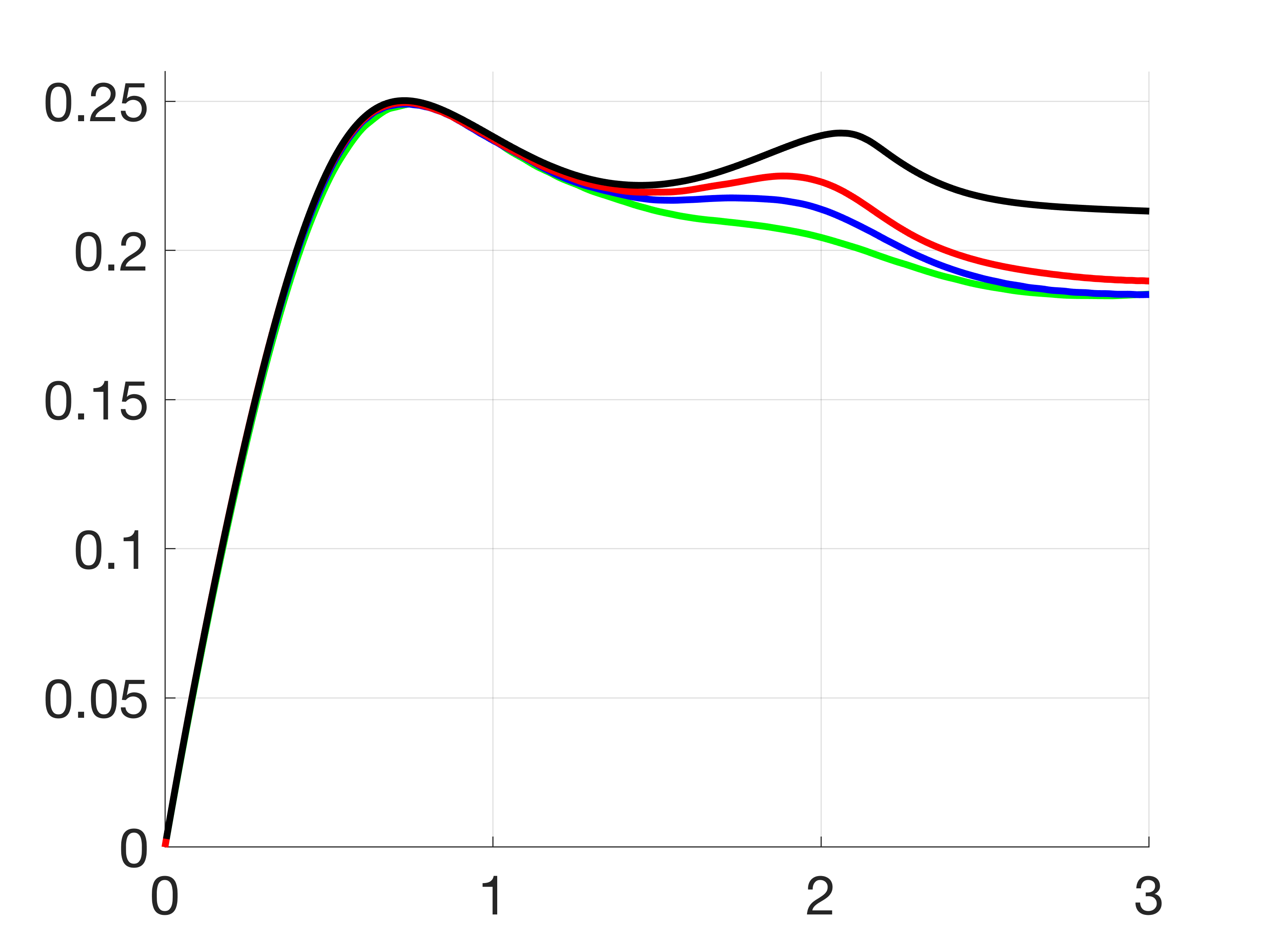}
  \caption{
    Circularity (left), center of mass (middle), and rising velocity
    (right) over time for test case 2.
  Green: $h=2^{-5}$. Blue: $h=2^{-6}$. Red: $h=2^{-7}$. Black:
  reference data from group 3 of \cite{Hysing09}.
(For interpretation of the colors in this figure, the reader is referred to the web version of this article.)}
\label{fig:bbt23}
\end{figure}

Furthermore, we plot the circularity, center of mass, and rise velocity over time 
in Figure~\ref{fig:bbt23}.
However, even for the finest grid $h=2^{-7}, \epsilon=0.005$, differences
remain visible. We conclude that the resolution is still not fine enough to get sufficiently
close to the reference solution.

\subsection{Rayleigh–Taylor instability}
The Rayleigh–Taylor instability is a two-phase instability 
which occurs whenever two fluids of different density 
are accelerated against each other. 
We consider a similar setting as in \cite{GuermondQuartapelle00}.
This problem consists of two layers of fluid initially at rest in the gravity field
in
the domain  $\Omega = [0,1/2]\times [-2,2]$.  The initial position of the
perturbed interface is 
$\eta(x) = -0.1\cos(2\pi x)$. The heavy fluid is above and the density
ratio is 3 ($\rho_1=3$, $\rho_2=1$). The (initial) transition between the two fluids
is regularized by hyperbolic tangent: 
\[
  \phi(t=0) = 
  \mathrm{tanh}\left(\frac{1}{\sqrt{2}\epsilon}(y+0.1\cos(2\pi x))\right).
\]
We take gravity $g=2$ so that the non-dimentional time scale is the same as 
the time scale of Tryggvason \cite{Tryggvason88}.
The viscosity in both fluids is taken to be $\nu = \sqrt{2}/Re$, where 
$Re$ is the Reyholds number. Both $Re=1000$, and $Re=5000$ are used in the
numerical simulations.
For the phase-field model parameters, we use $\epsilon=1.28 h$, $\gamma=10^{-3}
\epsilon$, and take a small surface tension constant $\sigma =
0.01\epsilon$. We notice that no surface tension effect is taken into account 
in the model \cite{GuermondQuartapelle00}.
The same flow boundary condition as the rising bubble problem is used
here, namely, 
the upper and lower boundaries are set to no-slip, 
and the left and right boundaries are set to slip conditions.
Again, we take polynomial degree $k=2$ and consider uniform rectangular
meshes with mesh size $h=2^{-6}$  and $h=2^{-7}$. 
%Due to the large Reynolds
%number for this problem, the flow is convection dominated.
We take variable time step size to be 
\[
  \delta t = 0.1 h/v_{\max},
\]
where $v_{\max}$ is an estimated maximum velocity magnitude at the current
time step.
We run the simulation till time $t=2.5$.
On the coarse mesh with $h=2^{-6}$,
a total of 1808 time steps is used 
%on the coarse mesh with $h=^{-6}$ when 
when
$Re=1000$, and 1966 steps when $Re=5000$.
On the fine mesh with $h=2^{-7}$, 
a total of 3688 time steps is used 
%on the fine mesh where $h=^{-7}$ when 
when
$Re=1000$, and 4240 steps when $Re=5000$.

\begin{figure}[h!]
\centering
\includegraphics[width=.12\textwidth]{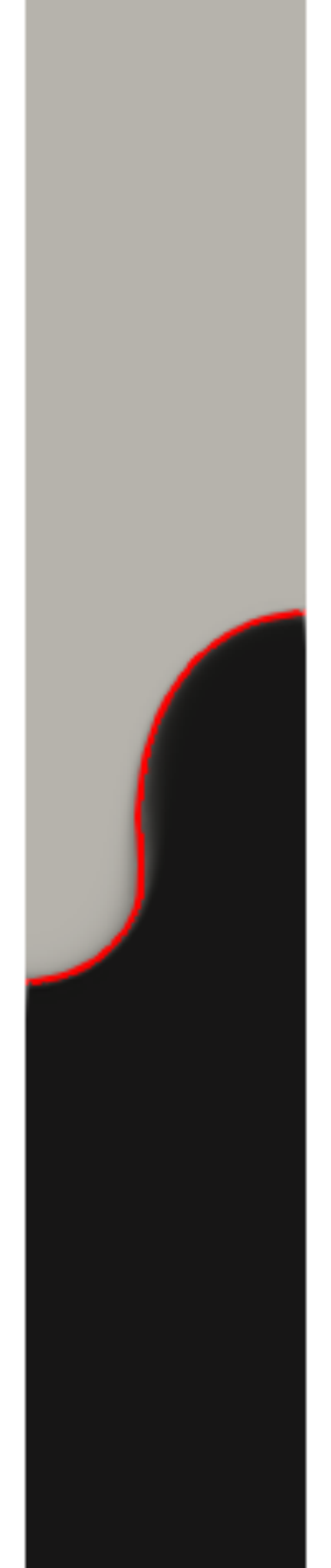}
\includegraphics[width=.12\textwidth]{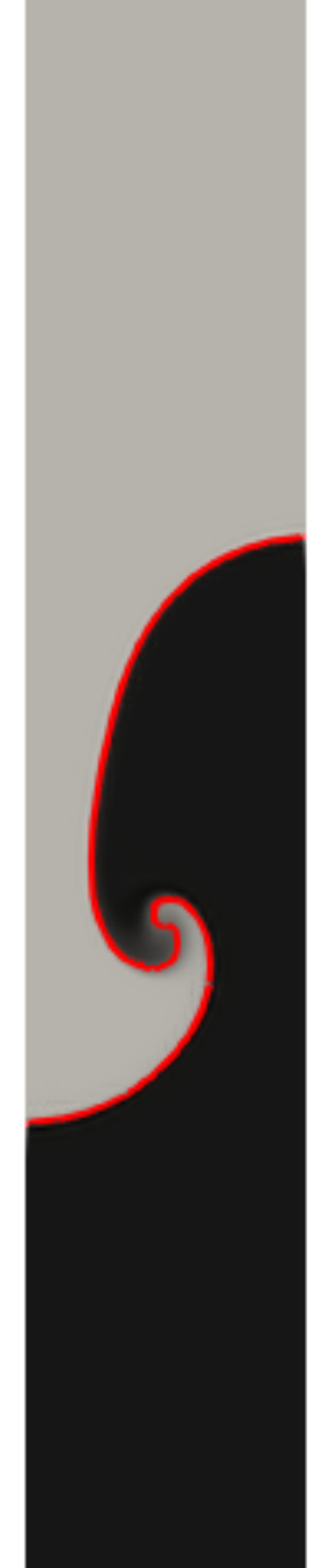}
\includegraphics[width=.12\textwidth]{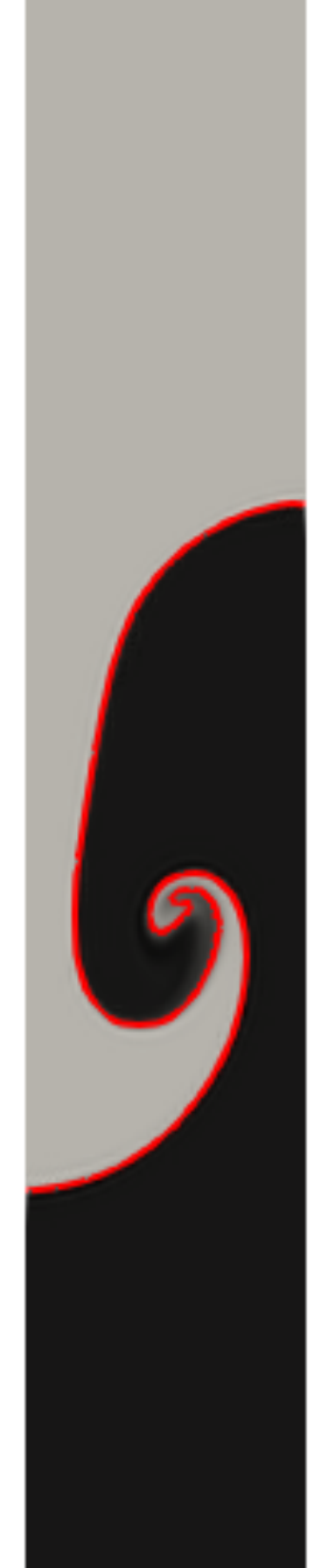}
\includegraphics[width=.12\textwidth]{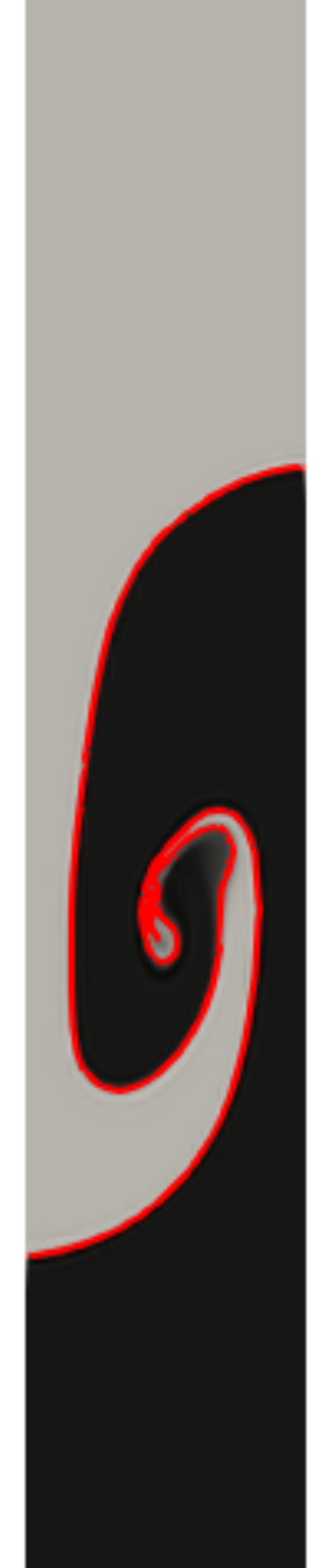}
\includegraphics[width=.12\textwidth]{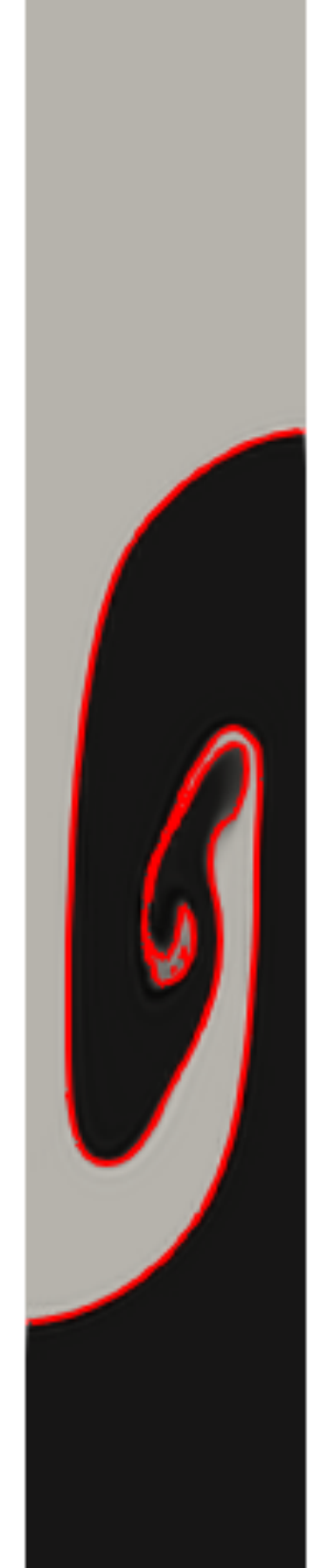}
\includegraphics[width=.12\textwidth]{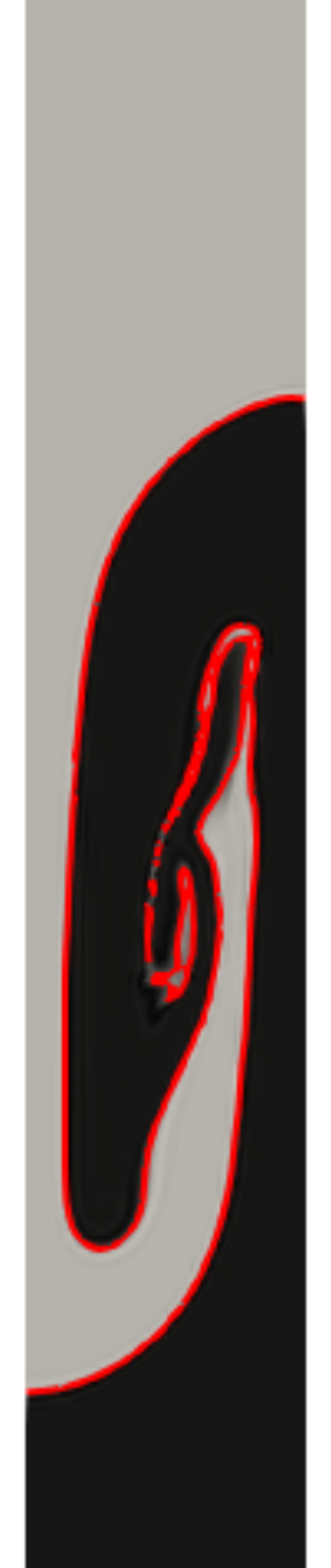}\\[.2ex]
\includegraphics[width=.12\textwidth]{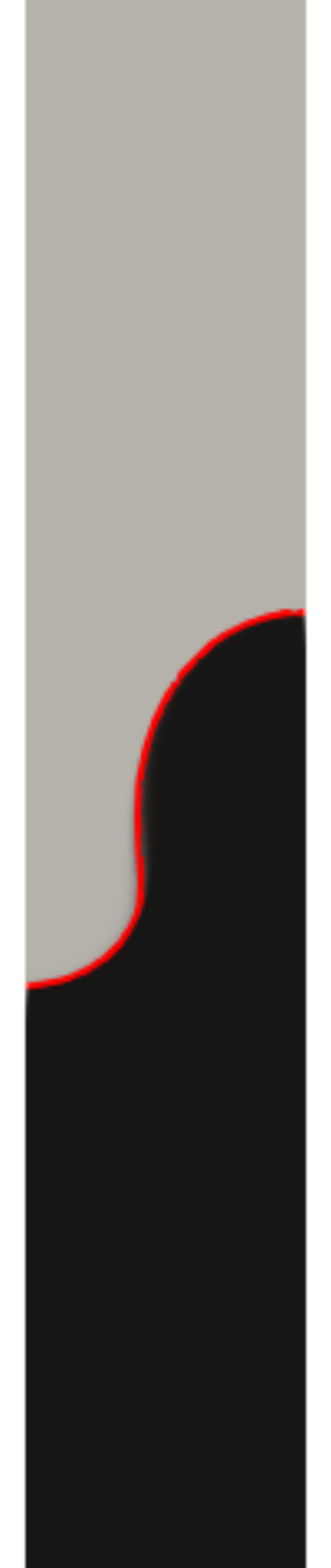}
\includegraphics[width=.12\textwidth]{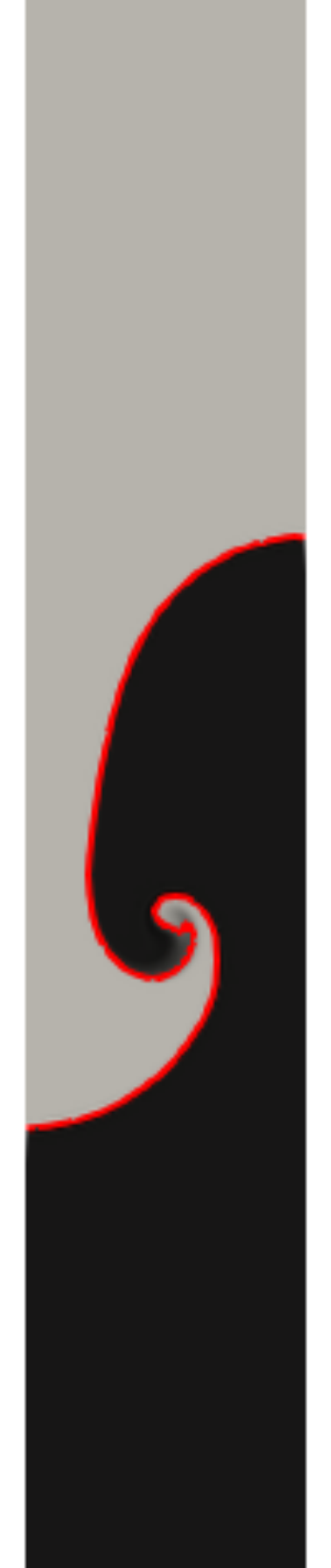}
\includegraphics[width=.12\textwidth]{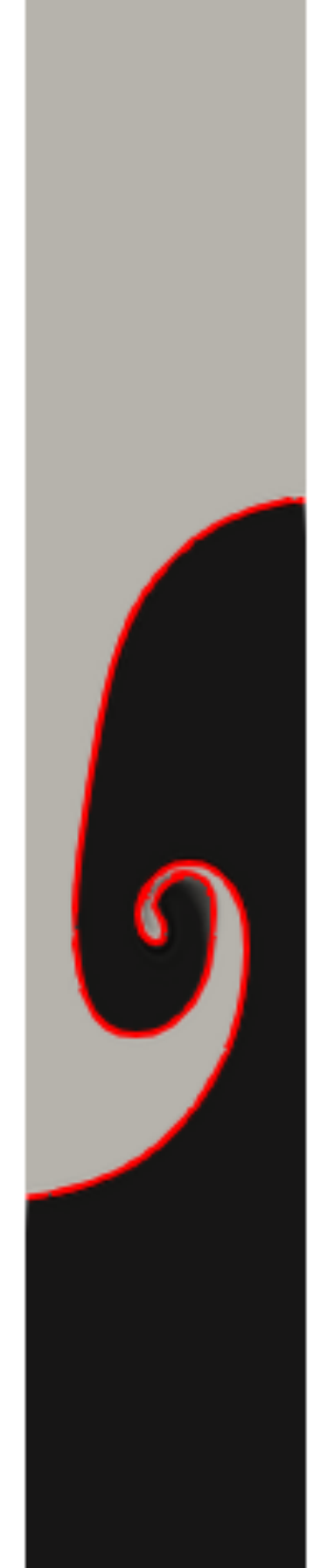}
\includegraphics[width=.12\textwidth]{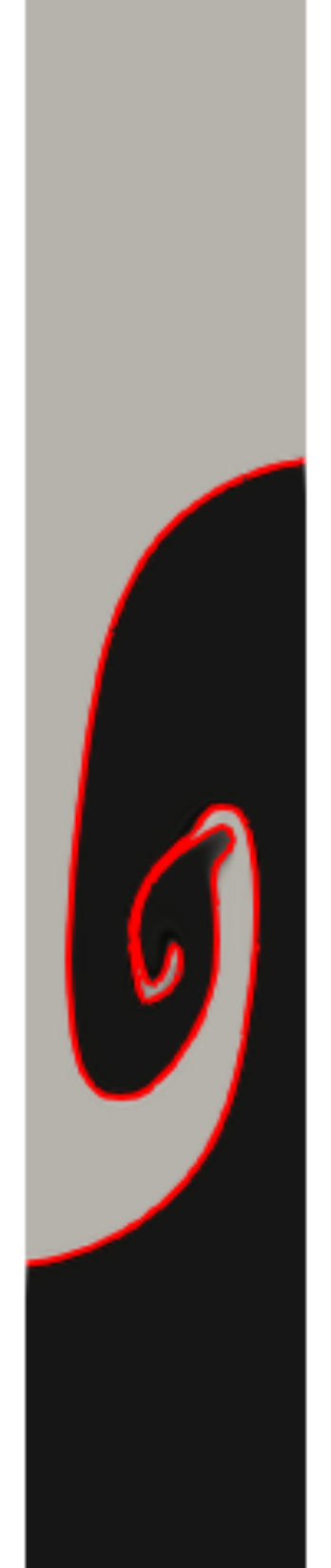}
\includegraphics[width=.12\textwidth]{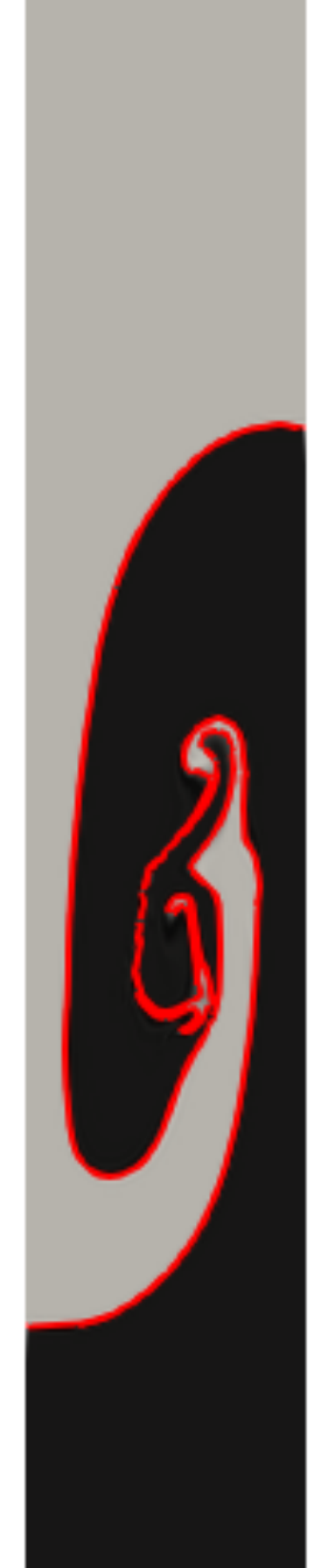}
\includegraphics[width=.12\textwidth]{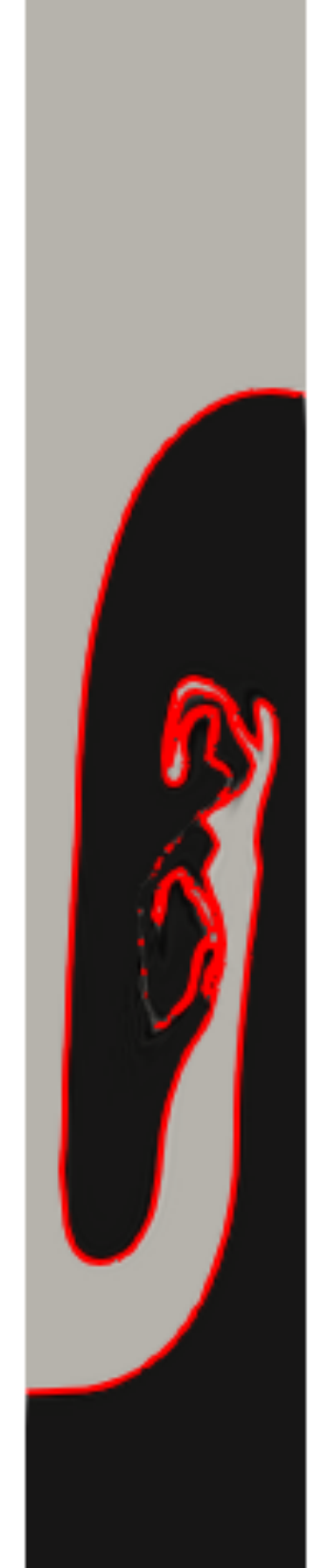}
\caption{
   $Re = 1000$. 
   Contour of the Rayleigh-Taylor instability problem
   on the subdomain $[0,0.5]\times [-1.41,1.41]$
   at time  $t=1,1.5,1.75,2.0,2.25,2.5$ (from left to right)
   Top: $h=2^{-6}$. Bottom: $h=2^{-7}$. 
  Red contour line: the interface $\phi_h=0$.
(For interpretation of the colors in this figure, the reader is referred to the web version of this article.)}
\label{fig:rt}
\end{figure}

\begin{figure}[h!]
\centering
\includegraphics[width=.12\textwidth]{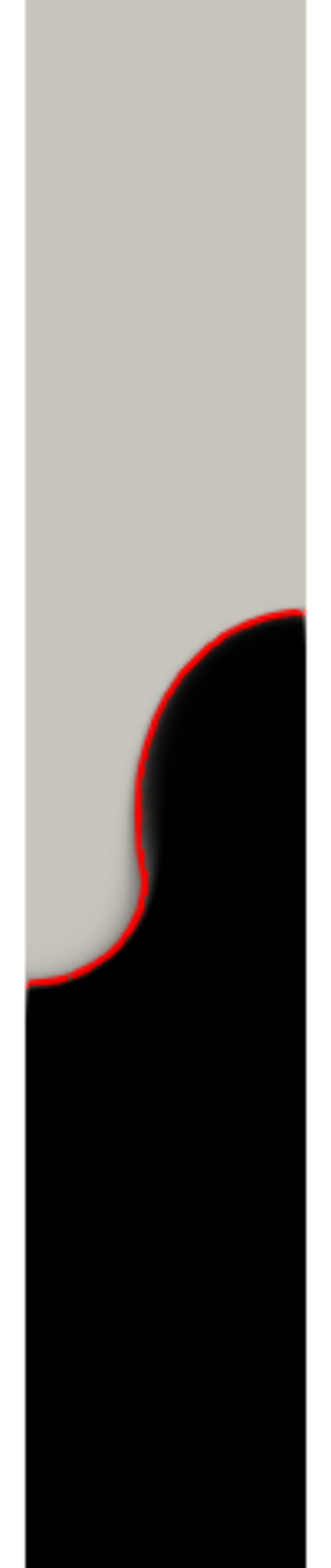}
\includegraphics[width=.12\textwidth]{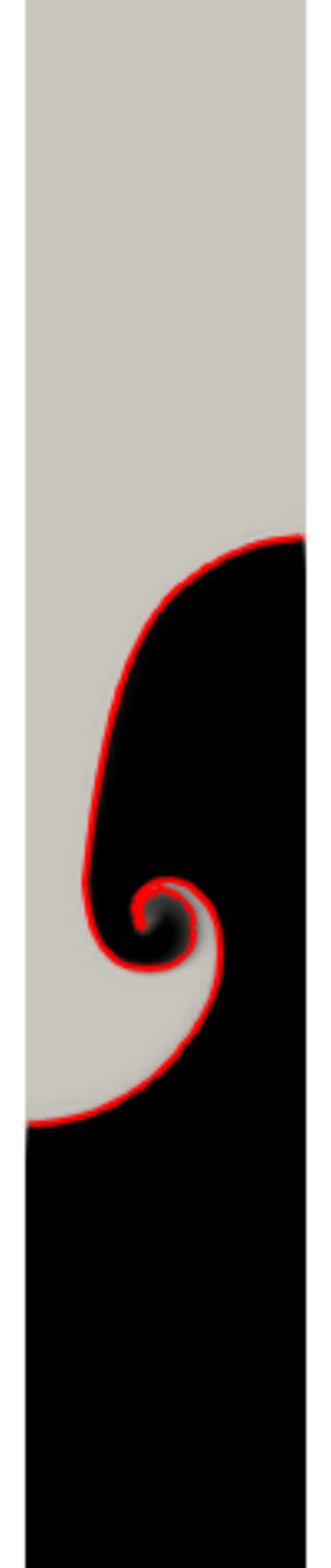}
\includegraphics[width=.12\textwidth]{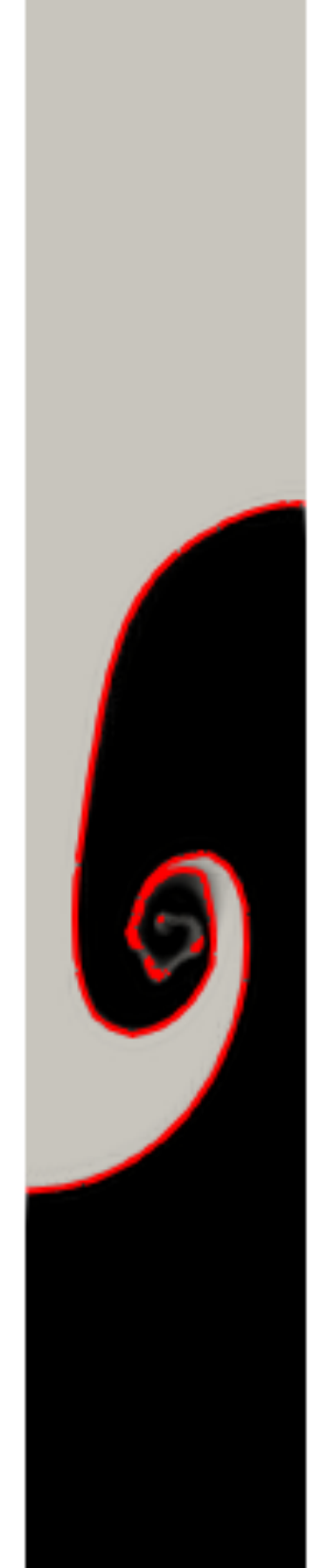}
\includegraphics[width=.12\textwidth]{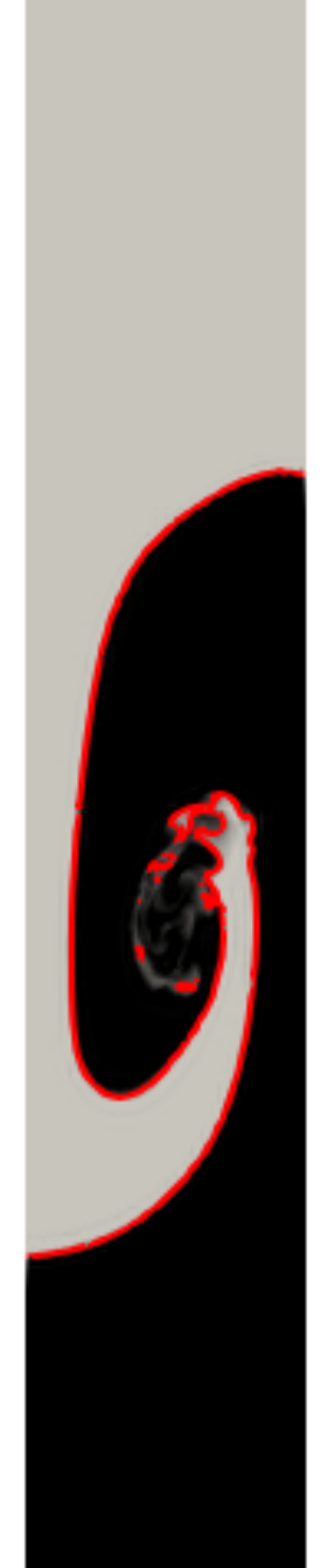}
\includegraphics[width=.12\textwidth]{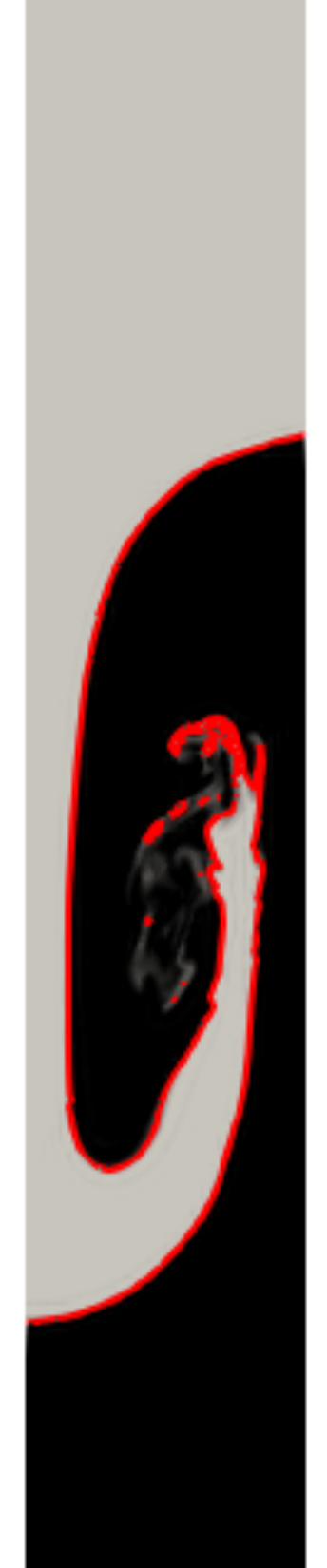}
\includegraphics[width=.12\textwidth]{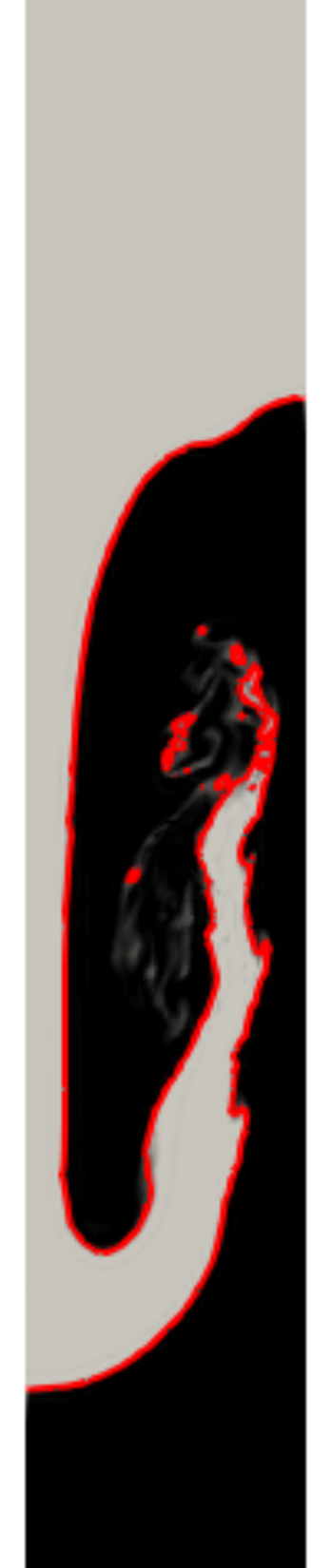}\\[.2ex]
\includegraphics[width=.12\textwidth]{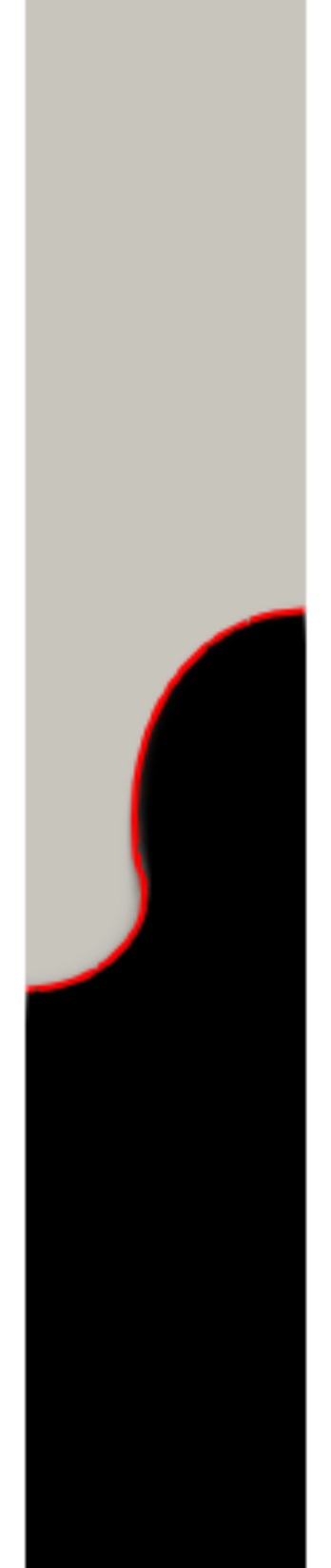}
\includegraphics[width=.12\textwidth]{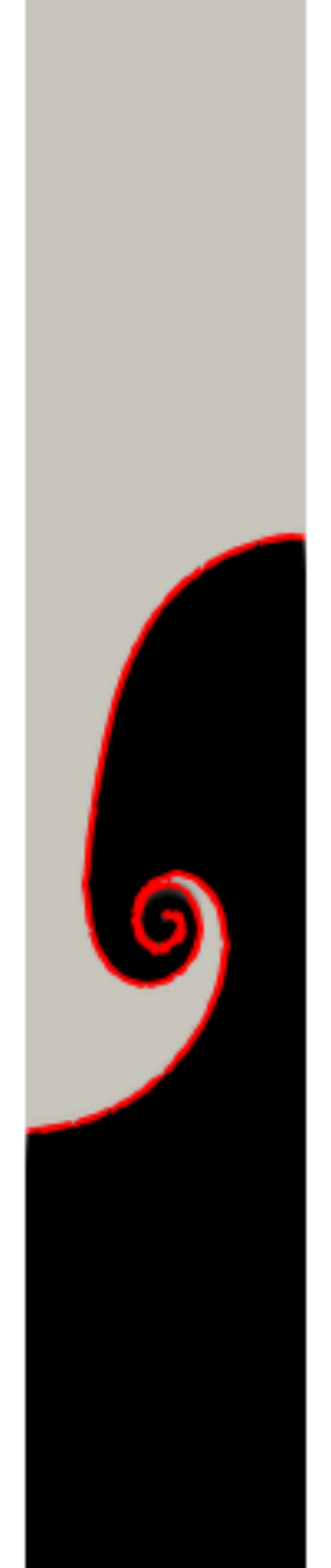}
\includegraphics[width=.12\textwidth]{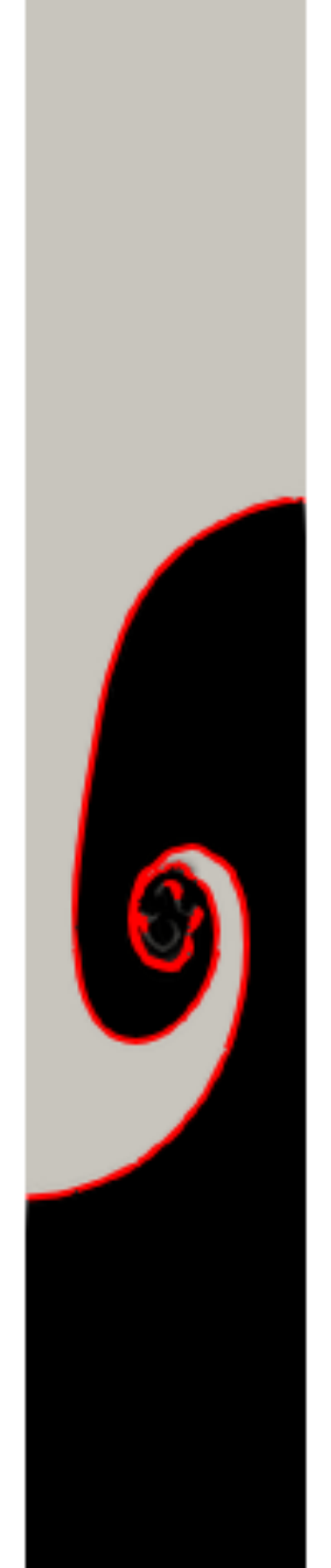}
\includegraphics[width=.12\textwidth]{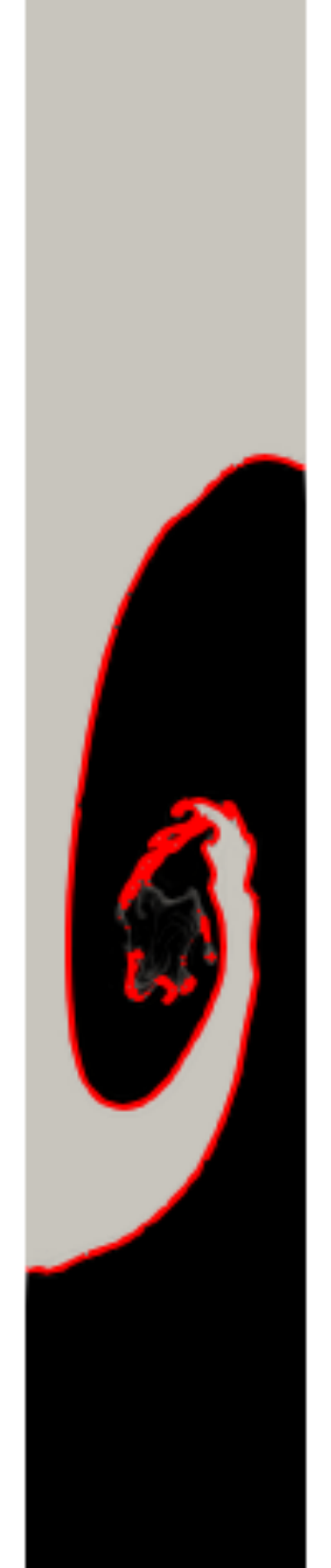}
\includegraphics[width=.12\textwidth]{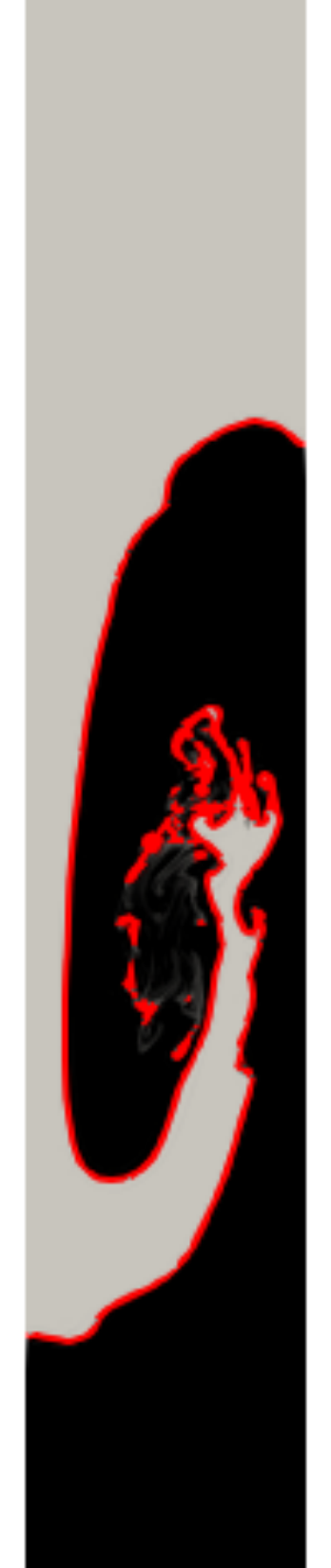}
\includegraphics[width=.12\textwidth]{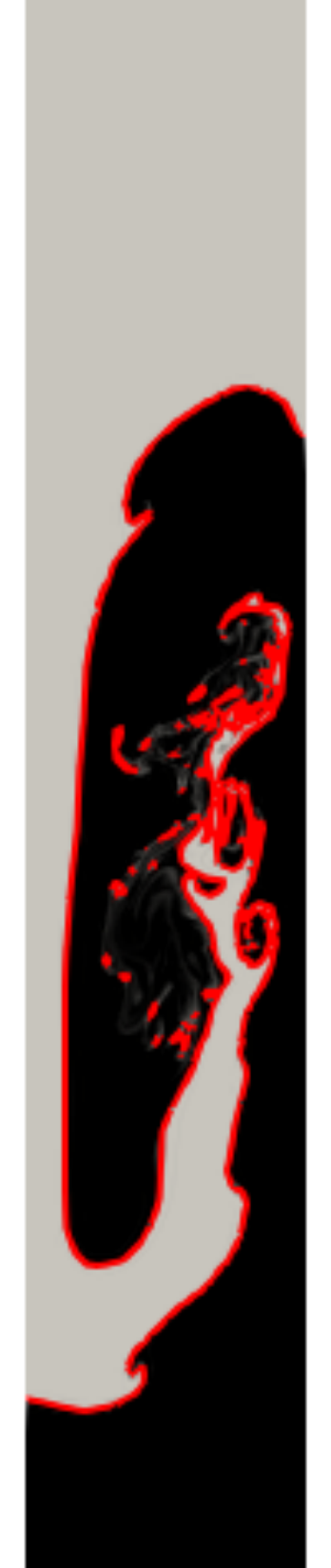}
\caption{
   $Re = 5000$. 
   Contour of the Rayleigh-Taylor instability problem
   on the subdomain $[0,0.5]\times [-1.41,1.41]$
   at time  $t=1,1.5,1.75,2.0,2.25,2.5$ (from left to right)
   Top: $h=2^{-6}$. Bottom: $h=2^{-7}$. 
   %Polynomial degree $k=2$.
  Red contour line: the interface $\phi_h=0$.
(For interpretation of the colors in this figure, the reader is referred to the web version of this article.)}
\label{fig:rt2}
\end{figure}

The time evolution of the contour of the phase-field variable $\phi_h$
for $Re = 1000$ and $Re=5000$ are plotted in Figure~\ref{fig:rt} and 
Figure~\ref{fig:rt2}, respectively, 
at times 1, 1.5, 1.75, 2, 2.25, 2.5.
From the results for $Re=1000$ in Figure~\ref{fig:rt}, we observe that 
the solutions on the two meshes are consistent in that they show similar
structures and differ only in fine details at large time.
Second, by comparing the solutions
for the results for $Re=5000$ in Figure~\ref{fig:rt2},  we observe 
that the solution on the two  meshes  
are in very good agreement in the early
stage of the time evolution $t \le 1.75$. 
Some noticeable differences occur at later times and
consist in the development of structures within the main 
vortex that are more complex on
the fine mesh than on the coarse one. 
%When comparing the solutions at $t=2.5$, we observe
%some disagreement in the shape of the roll-up of the rising bubble, but fair agreement in
%the overall shape (external and internal) of the falling bubble.
All these results are qualitatively similar to the results in 
\cite{GuermondQuartapelle00}, where a projection FEM was used to solve 
the variable density flow without surface tension.

Finally in Table \ref{tab:rt}, we compare the minimal and maximum of $y$ 
position (bubble and spike location, respectively) of the interface $\phi_h=0$ 
with results in
\cite{GuermondQuartapelle00} (data in \cite{GuermondQuartapelle00} is 
extrapolated up to two digits accuracy from Figure~4 therein).
Excellent agreement of our results with the reference data is observed.
In particular, we observe that, on the same mesh, the bubble/spike
locations are very similar for $Re=1000$ and $Re=5000$.
\begin{table}[ht!]
\begin{center}
  %\footnotesize
  %\scalebox{1}
  {
  \begin{tabular}{c c  c  c  c c c c c} \hline
    &  $h$ & $Re$ & t=1 &t=1.5& t=1.75&t=2.0&t=2.25&t=2.5\\
\hline
    & $2^{-6}$ & 1000 & -0.3617&-0.6139&-0.7351&-0.8511&-0.9706&-1.0970\\ 
    & $2^{-7}$ & 1000 &
    -0.3690&-0.6235&-0.7460&-0.8654&-0.9763&-1.0957\\[.5ex]
bubble
%   $y|_{\phi_h(y,0)=0} 
    & $2^{-6}$ & 5000 & 
     -0.3617 &  -0.6142  & -0.7358&   -0.8513 &  -0.9712 &  -1.0974\\
    & $2^{-7}$ & 5000 & 
       -0.3690 &  -0.6236 &  -0.7461 &  -0.8655&   -0.9781 &  -1.0963\\
    \hline
    & ref      & 1000 &
    -0.37&-0.62&-0.74&-0.86&-0.98&-1.11\\
    & ref      & 5000 &
    -0.39&-0.64&-0.75&-0.87&-0.98&-1.11\\
    \hline 
    \hline
    & $2^{-6}$ & 1000 &
     0.2959  &  0.4312  &  0.4946 &   0.5577  &  0.6203    &0.6836\\
    & $2^{-7}$ & 1000 &
    0.2988 &   0.4347   & 0.5009  &  0.5671  &  0.6301  &  0.6880\\[.5ex]
spike
%   $y|_{\phi_h(y,0)=0} 
    & $2^{-6}$ & 5000 &
    0.2959  &  0.4314  &  0.4949 &   0.5581  &  0.6207  &  0.6841\\
    & $2^{-7}$ & 5000 &
        0.2988  &  0.4347  &  0.5010  &  0.5672  &  0.6323  &  0.6933\\
    \hline
    & ref      & 1000 & 
        0.30  &  0.42  &  0.48 &   0.55 &   0.61&    0.69\\
    & ref      & 5000 &
        0.30 &   0.43  &  0.50 &   0.58  &  0.65 &   0.72\\
    \hline 
\end{tabular}
}
\end{center}
\caption{\it 
  Bubble and spike locations for the Rayleigh-Taylor instability problem 
at different time.}
\label{table:rt}
%\end{wraptable}
\end{table}

\section{Conclusion and future work}
\label{sec:conclude}
We have presented a  novel divergence-free 
HDG scheme for 
a Cahn-Hilliard phase-field model for two-phase incompressible flow.
The (linear and decoupled) 
fully discrete is observed to be second-order accurate in time, and optimal
order accurate in space ($(k+1)$-th order for the $L^2$-errors when
polynomials of degree $k$ is used).
Benchmark results are presented
for the classical
bubble-rising problem and the Rayleigh-Taylor instability problem, 
which are consistent with results in the
literature.

This work consists of our initial investigation of divergence-free HDG
schemes for phase-field model of two-phase incompressible flow. 
The scheme is robust in the convection-dominated regime, produce a globally
divergence-free velocity approximation, and can be efficiently 
implemented via static condensation.
In the
future, we plan to investigate on the derivation and analysis of efficient 
fully discrete energy stable schemes for phase-field model of two-phase
incompressible flow, where we shall consider models admit energy law such
as
those in \cite{ShenYang10,Abels12}.  
The combination of our (spatial) divergence-free HDG scheme with the 
recently introduced scalar auxilary variable (SAV) technique \cite{SXY19} for the phase-field
temporal discretization seems to be a very promising approach, which
consists of our ongoing work.
%Detailed convergence 
%analysis of the scheme consists of our ongoing work.
%We are also planning to investigate 
%the extension of this scheme to fluid-structure interactions.
\bibliographystyle{siam}

\end{document}